\documentclass[12pt]{amsart}
\usepackage{amssymb}
\ifx\mathrm\undefined\let\mathrm\rm\fi
\ifx\mathbf\undefined\let\mathbf\bf\fi
\ifx\mathfrak\undefined\let\mathfrak\frak\fi
\ifx\mathcal\undefined\let\mathcal\cal\fi
\ifx\mathbb\undefined\let\mathbb\Bbb\fi
\ifx\emph\undefined\let\emph\it\fi

\newcommand{\Z}{{\mathbb Z}}

\newcommand{\R}{{\mathbb R}}
\newcommand{\C}{{\mathbb C}}
\newcommand{\Q}{{\mathbb Q}}

\newcommand{\Ref}[1]{{(\ref{#1})}}
\newcommand{\be}{\begin{displaymath}}
\newcommand{\ee}{\end{displaymath}}
\newcommand{\bea}{\begin{eqnarray*}}
\newcommand{\eea}{\end{eqnarray*}}

\font\bb=msbm10 at9.98pt
\def\semidirect{\hbox{$\;$\bb\char'156$\;$}}

\newcommand{\dontprint}[1]{\relax}

\newenvironment{prf}{\noindent{\it Proof\/}:}{$\;\square$
\par\medskip}

\newtheorem%
{thm}{Theorem}[section]
\newtheorem%
{proposition}[thm]{Proposition}
\newtheorem%
{lemma}[thm]{Lemma}

\newtheorem%
{lemmadef}[thm]{Lemma-Definition}
\newtheorem%
{corollary}[thm]{Corollary}
\newtheorem%
{conjecture}[thm]{Conjecture}

\title[The elliptic gamma function]{The elliptic gamma function and $\mathrm{SL}(3,\Z)\semidirect\Z^3$}
\author[G. Felder and A. Varchenko]
{Giovanni Felder${}^{*}$ 
\and Alexander Varchenko${}^{**,1}$}
\thanks{${}^1$Supported in part by NSF grant  DMS-9801582}
\date{June 1999}
\begin{document}
\maketitle
\centerline{\it ${}^*$Departement Mathematik, ETH-Zentrum,}
\centerline{\it 8092
Z\"urich, Switzerland}
\centerline{felder@math.ethz.ch}
\medskip
\centerline{\it ${}^{**}$Department of Mathematics,
University of North Carolina at Chapel Hill,}
\centerline{\it Chapel Hill, NC 27599-3250, USA}
\centerline{av@math.unc.edu}
\newcommand{\sig}{\sigma}
\begin{abstract}
The elliptic gamma function is a  generalization
of the Euler gamma function and is  associated to an elliptic curve.
Its trigonometric and rational degenerations are the Jackson
q-gamma function and the Euler gamma function, respectively.
The elliptic gamma function appears in Baxter's formula
for the free energy of the eight-vertex model and
in the hypergeometric solutions of the elliptic qKZB equations.
In this paper, the properties of this function are studied.
In particular we show that elliptic gamma functions are 
generalizations of automorphic forms of $G=SL(3,\Z)\semidirect\Z^3$ 
associated to a non-trivial class in $H^3(G,\Z)$.
\end{abstract}
\section{Introduction}
This paper
deals with the properties of the {\em elliptic
gamma function}, an elliptic generalization
of the Euler gamma function. It is the meromorphic function
of three complex variables $z,\tau,\sig$, with 
$\mathrm{Im}\,\tau,\mathrm{Im}\,\sig>0$ defined
by the convergent infinite product
\[
\Gamma(z,\tau,\sig)=\prod_{j,k=0}^\infty
\frac
{1-e^{2\pi i((j+1)\tau+(k+1)\sig-z)}}
{1-e^{2\pi i(j\tau+k\sig+z)}}\,.
\]
It is the unique solution
of a functional equation involving a Jacobi
theta function: 
Let for $z,\tau\in\C$ with
$\mathrm{Im}\,\tau>0$, $\theta_0$ denote
the theta function
\[
\theta_0(z,\tau)=\prod_{j=0}^\infty
(1-e^{2\pi i((j+1)\tau -z)}) 
(1-e^{2\pi i(j\tau +z)}).
\]
Then the elliptic gamma function may be characterized as
follows.

\begin{thm}\label{t-1}
Suppose that $\tau,\sig$ are complex numbers with positive 
imaginary part. Then $u(z)=\Gamma(z,\tau,\sig)$ is the
unique meromorphic solution of the difference equation
\[
u(z+\sig)=\theta_0(z,\tau)u(z)
\]
such that:
\begin{enumerate}
\item[(i)] $u(z)$ obeys $u(z+1)=u(z)$ and is holomorphic
on the upper half plane
$\mathrm{Im}\, z>0$,
\item[(ii)] $u((\tau+\sig)/2)=1$. 
\end{enumerate}
\end{thm}

This Theorem is proved in \ref{s-egf} below.

The elliptic gamma function was introduced by Ruijsenaars \cite{R}.
 Similar double products appeared earlier in statistical
mechanics.
Probably the first (implicit) appearance 
of the elliptic gamma function in this context
is in Baxter's formula \cite{B} for the
free energy of the eight-vertex model. This model
has four parameters, which can be taken to be $c_B,x_B,q_B,z_B$
in Baxter's notation (to which we add a subscript B to avoid conflicts 
with our  notations), see \cite{B} eqs.\ (D1)-(D8). If we set $c_B=c$,
$q_B=e^{2\pi i\tau}$, $x_B=e^{2\pi i\sig}$, $z_B=e^{2\pi iu}$, then
 the free energy (times the
inverse temperature) is
$ f(c,u,\tau,\sig)=-\ln\,c-\ln\,Z(u,\tau,\sig)$ where $Z$ can be
expressed in terms of theta and elliptic gamma functions:
\[
Z(u,\tau,\sig)
=
\frac
{\theta_0(\sig-u,2\tau)\theta_0(\sig+u,2\tau)
\Gamma(\sig-u,\tau,4\sig)\Gamma(\sig+u,\tau,4\sig)}
{\theta_0(\tau,2\tau)\theta_0(2\sig,2\tau)
\Gamma(3\sig-u,\tau,4\sig)\Gamma(3\sig+u,\tau,4\sig)}
\,.\]
The elliptic 
gamma function and similar double and even triple infinite products appear in
correlation functions of the eight-vertex model \cite{Kyoto1, Kyoto2}
and boundary spontaneous magnetization at corners of the 
Ising model \cite{DP}.

Our own interest in the elliptic gamma functions arose from
the study of
hypergeometric solutions of elliptic qKZB difference equations
\cite{FTV,FV}.
In these solutions the role of powers of linear functions appearing
in the Gauss hypergeometric function is played by ratios
of elliptic gamma functions.

In this paper, after reviewing 
some well-known properties
of the theta function $\theta_0$,
 we derive several identities for elliptic gamma functions.
The most remarkable identities
are the ``modular'' three-term relations of Theorem \ref{t-4},
which connect values of $\Gamma$ at points related by SL$(3,\Z)$ acting
on the periods $(\tau,\sig)$ by fractional linear transformations.
For example, there is a polynomial $Q(z;\tau,\sig)$ of degree three
in $z$ whose coefficients are rational functions of $\tau,\sig$ such
that
\[
\Gamma(
 z/\sig,\tau/\sig,-1/\sig)=
e^{i\pi Q(z;\tau,\sig)}
{\Gamma(
{({z-\sig})/\tau,-1/\tau,-\sig/\tau})}
{\Gamma(z,\tau,\sig)}.
\]
These identities have an interpretation in terms of a generalization
of Jacobi modular forms:
$\Gamma$ may be interpreted as the value of generators
of an ``automorphic form of
degree 1'': just as the theta function $\theta_0$ is a ``degree 0''
automorphic form  associated to a 1-cocycle in $H^1(G,M)$
where $G=SL(2,\Z)\semidirect\Z^2$ and $M=\exp 2\pi i\Q(\tau)[z]$, the
elliptic gamma function defines a ``degree 1'' automorphic form
associated to a $2$-cocycle in $H^2(G,M)$
with $G=\mathrm{SL}(3,\Z)\semidirect\Z^3$ and $M=\exp 2\pi i\Q(\tau,\sigma)[z]$.
See Sect.\ \ref{s-last} for a more precise statement.

The modular identities have interesting degenerations:
limiting versions relate the gamma functions
at points where the periods are linearly dependent over the rationals
to the Euler dilogarithm function (Theorems \ref{t-5}, \ref{t-6}, \ref{formula}). 
On the other hand, in the
semiclassical limit $\epsilon\to 0$ 
 the ``phase function'' $\Gamma(z+\epsilon,\tau,2\epsilon/\beta)/
\Gamma(z-\epsilon,\tau,2\epsilon/\beta)$ tends to $\theta_0(z,\tau)^\beta$,
and the modular identities reduce to the Jacobi modular properties of 
theta functions, see Section \ref{S-6}.

One intriguing aspect of the elliptic gamma function, which is
a priori defined for periods $\tau,\sig$ in the upper half plane,
is that it may be extended, by a simple reflection, to a function
of $\tau,\sig\in\C-\R$ in such a way that all identities remain
true. Moreover the elliptic gamma function
 also has a limit as $\tau$ or $\sig$ (but not
both) approach, from either side,
a subset $X$  of full Lebesgue measure 
of the real axis, see Theorem \ref{t-real}. 
This subset contains all irrational algebraic real
numbers.

The paper is organized as follows. In Sect.\ \ref{theta} we
review some well-known properties of the odd Jacobi theta function.
Then in Sect.\ \ref{S-3} we introduce the elliptic gamma
function: after giving its definition and its elementary
properties, we study its trigonometric and rational degenerations.
Then we derive a summation formula that allows us to study
the limit as one approaches the real axis.
The modular properties of the elliptic gamma function are
given in Sect.\ \ref{S-4}. The properties of the elliptic
gamma function at special values of its arguments are
studied in Sect.\ \ref{S-5}. There the relation to
dilogarithms appears. In Sect.\ \ref{S-6} we
study the semiclassical limit, in which the identities
of elliptic gamma functions reduce to differential and
difference equations obeyed by powers of theta functions.
In the last section we introduce the notion of automorphic forms
of degree 1 and 
relate  the modular identities satisfied by gamma functions to
these 
automorphic forms and to the cohomology of $\mathrm{SL}(3,\Z)\semidirect\Z^3$.

\section{Theta functions}\label{theta}
\subsection{The theta function}\label{A-theta}
Jacobi's first theta function is defined by the series
\[
\theta(z,\tau)=-\sum_{j\in\Z}
e^{i\pi\tau(j+1/2)^2+2\pi i(j+1/2)(z+1/2)},\qquad z,\tau\in \C,\qquad 
\mathrm{Im}\,{}\tau>0.
\]
It is an entire holomorphic odd function such that
\begin{equation}\label{eq-frt}
\theta(z+n+m\tau,\tau)=(-1)^{m+n}e^{-\pi im^2\tau-2\pi i mz}\theta(z,\tau),
\qquad m,n\in\Z,
\end{equation}
and obeys the heat equation 
\begin{equation}\label{e-heat}
4\pi i\textstyle{\frac\partial{\partial\tau}}\theta(z,\tau)=
\theta''(z,\tau).
\end{equation}
Its transformation properties with respect to $SL(2,\Z)$ are described
in terms of generators by the identities:
\[
\theta(-z,\tau)=-\theta(z,\tau),\qquad
\theta(z,\tau+1)=e^{\frac{i\pi}{4}}\theta(z,\tau),
\qquad
\theta\left(\frac{z}{\tau},
-\frac{1}{\tau}
\right)=-i
{\sqrt{-i\tau}}
e^{\frac{i\pi z^2}{\tau}}
\theta(z,\tau).
\]
The square root is the one in the right
half plane.

\subsection{Infinite products}
Let $x,q\in\C$  with $|q|<1$. The function
\[
(x;q)=\prod_{j=0}^\infty(1-xq^j)
\]
is a solution of the functional equation
\[
(qx;q)=\frac 1{1-x}(x;q).
\]
Using the identity 
\begin{equation}\label{e-Id}
1-y=\exp(-\sum_{j=1}^\infty y^j/j), \qquad |y|<1,
\end{equation}
and
summing the geometric series, yields the summation formula
\begin{equation}\label{e-tu}
(x;q)=\exp\left(-\sum_{j=1}^\infty\frac{x^j}{j(1-q^j)}\right),
\qquad |x|<1,\quad |q|<1.
\end{equation}

\subsection{Product representation of theta functions}
Let  $x=e^{2\pi i z}$ and  $q=e^{2\pi i \tau}$. Then
we have the Jacobi triple product identity
\begin{equation}\label{eq-thetafunction}
\theta(z,\tau)=ie^{\pi i(\tau/4-z)}(x;q)(q/x;q)(q;q)
\end{equation}
In particular we have $\theta'(0,\tau)=2\pi\,\eta(\tau)^3$,
with  $\eta(\tau)=e^{\pi i\tau/12}(q;q)$ the Dedekind function.

{}From the modular properties of $\theta$, we deduce the 
modular properties
of the Dedekind function:
$\eta(\tau+1)=e^{\pi i/12}\eta(\tau)$ and $\eta(-1/\tau)=
(-i\tau)^{1/2}\eta(\tau)$ (up to a third root of unity, which
is $1$, as one sees by setting $\tau=i$).

We will also need
\[
\theta_0(z,\tau)=(x;q)(q/x;q)=-i\frac{e^{\pi i(z-\tau/4)}}{(q;q)}\theta(z,\tau).\]
This function obeys
\begin{eqnarray}\label{EGF1}
\theta_0(z+1,\tau)&=& \theta_0(z,\tau),\nonumber
\\
\theta_0(z+\tau,\tau)&=& -e^{-2\pi iz}\theta_0(z,\tau),
\\
\theta_0(\tau-z,\tau)&=&\theta_0(z,\tau).\nonumber
\end{eqnarray}
Its modular properties follow from those of $\theta$ and $\eta$: they are
$\theta_0(z,\tau+1)=\theta_0(z,\tau)$, 
and if $z'=z/\tau$, $\tau'=-1/\tau$,
\begin{equation}\label{e-modt0}
e^{\pi i(\tau/6-z)}\theta_0(z,\tau)=
i\, e^{\pi i(-zz'+\tau'/6-z')}\theta_0(z',\tau').
\end{equation}
The summation formula \Ref{e-tu} implies
\begin{equation}\label{e-tsf}
\theta_0(z,\tau)=\exp\left(
-i\sum_{j=1}^\infty 
\frac{\cos\,\pi j(2z-\tau)}
{j\sin\,\pi j\tau}
\right),\qquad 0<\mathrm{Im}\,z<\mathrm{Im}\,\tau.
\end{equation}
\section{Elliptic gamma functions}\label{S-3}
\subsection{Definitions and elementary properties}\label{s-egf}
Here we consider two parameters $\tau$ and $\sig$ in the upper half
plane, and set $q=e^{2\pi i \tau}$, $r=e^{2\pi i \sig}$, and consider
the function of $x=e^{2\pi i z}$, 
\[
(x;q,r)=\prod_{j,k=0}^\infty(1-xq^jr^k)=(x;r,q)
\] It is a solution of the functional equations
\begin{equation}\label{e-fe}
(qx;q,r)=\frac{(x;q,r)}{(x;r)},
\qquad
(rx;q,r)=\frac{(x;q,r)}{(x;q)}.
\end{equation}
For $|x|<1$, we have the formula
\begin{equation}\label{eq-SumFormula}
(x;q,r)=\exp\left(-\sum_{l=1}^\infty x^l/l(1-q^l)(1-r^l)\right).
\end{equation}
It is obtained as \Ref{e-tu} 
by expanding the logarithm in a Taylor series
in $x$ and then summing the resulting geometric series.

The {\em elliptic gamma function} is 
\[
\Gamma(z,\tau,\sig)=
\frac{(qr/x;q,r)}{(x;q,r)}\,.
\]
\begin{thm}\label{t-2}
The elliptic gamma function  obeys the identities
\begin{eqnarray}\label{eq-Gamma9}
  \Gamma(z,\tau,\sig)&=&\Gamma(z,\sig,\tau),\\
  \Gamma(z+1,\tau,\sig)&=&\Gamma(z,\tau,\sig),\\
  \Gamma(z+\sig,\tau,\sig)&=&
  \theta_0(z,\tau)\Gamma(z,\tau,\sig),\\
  \Gamma(z+\tau,\tau,\sig)&=&
  \theta_0(z,\sig)\Gamma(z,\tau,\sig),
\end{eqnarray}
 and is normalized by
$\Gamma((\tau+\sig)/2,\tau,\sig)=1$.
As a function of $z$, $\Gamma(z,\tau,\sig)$ 
is a meromorphic function whose zeros and poles are all
simple. The zeros are at
$
z=(j+1)\tau+(k+1)\sig+l,
$
and the poles are at $z=-j\tau-k\sig+l$, where $j,k$ run
over nonnegative integers and $l$ over all integers.
\end{thm}

\medskip

\begin{prf}
It is obvious that $\Gamma$ is symmetric under interchange of
$\tau$ and $\sig$ and is 1-periodic. 
The remaining identities follow from
\Ref{e-fe}. 
The zeros of $(x;q,r)$ are at $x=q^{-j}r^{-k}$, $j,k=0,1,2,\dots$.
This implies the statement about zeros and poles.
\end{prf}

\medskip
\noindent{\it Proof of Theorem \ref{t-1}:} It remains to prove uniqueness.
The point is that $u(z)=\Gamma(z,\tau,\sig)$ has no zeros 
in the strip $0<\mathrm{Im}\,z<\mathrm{Im}\,\sig+\epsilon$, for some
$\epsilon>0$.
If $v(z)$ is another 1-periodic solution, holomorphic in the upper half
plane, then $v(z)/u(z)$ is a doubly periodic function with periods
$1$ and $\sig$. It  is holomorphic in the same strip and thus, by
periodicity, a bounded entire function. By Liouville's theorem,
$v/u$ is thus constant, which implies our claim. $\square$

\medskip
Finally, we mention some elementary identities, which may be thought of as the
analogues of the classical
formula $\Gamma(1-z)\Gamma(z)=\pi/\sin(\pi z)$:
\begin{proposition}\label{p-3}
\[
\Gamma(z,\tau,\sig)\Gamma(\sig-z,\tau,\sig)=\frac1{\theta_0(z,\sig)},
\qquad
\Gamma(z,\tau,\sig)\Gamma(\tau-z,\tau,\sig)=\frac1{\theta_0(z,\tau)},
\]\[
\Gamma(z,\tau,\sig)\Gamma(\tau+\sig-z,\tau,\sig)=1.
\]
\end{proposition}

\subsection{Trigonometric and rational limit}
We have the trigonometric and rational limit of the
theta function:
\[
\frac{\theta_0(\sig s,\tau)}{\theta_0(\sig,\tau)}
\stackrel{\tau\to i\infty}{\longrightarrow}
\frac
{1-e^{2\pi i\sig\, s}}
{1-e^{2\pi i\sig}}
\stackrel{\sig\to 0}{\longrightarrow}
s.
\]
Let
$\bar\Gamma$ be the function
\[
\bar\Gamma(s,\tau,\sig)=\frac{(r;r)}{(q;q)}
\theta_0(\sig,\tau)^{1-s}
\Gamma(\sig s,\tau,\sig),\qquad q=e^{2\pi i\tau},
r=e^{2\pi i\sig}.
\]
Then $u(s)=\bar\Gamma(s,\tau,\sig)$ is a solution of the functional equation
\[
u(s+1)=\frac{\theta_0(\sig s,\tau)}{\theta_0(\sig,\tau)}\,u(s).
\]
The normalization was chosen here so that $u(1)=1$.
As $\tau\to i\infty$ we recover F. H. Jackson's  ``q-gamma function'',
\[
\Gamma_{\mathrm{trig}}(s,\sig)=
\lim_{\tau\to i\infty}\bar\Gamma(s,\tau,\sig)=
(1-r)^{1-s}\frac{(r;r)}{(r^s;r)}\, .
\]
This function obeys the functional equation
$\Gamma_\mathrm{trig}(s+1,\sig)=\frac
{1-e^{2\pi i\sig\, s}}
{1-e^{2\pi i\sig}}
\Gamma_{\mathrm{trig}}(s,\sig)$,
and degenerates to
the Euler gamma function $\Gamma_\mathrm{Euler}(s)=\int_0^\infty
t^{s-1}e^{-t}dt$:
\[
\lim_{\sig\to 0}\Gamma_{\mathrm{trig}}(s,\sig)=\Gamma_{\mathrm{Euler}}(s).
\]
See \cite{A} for an account of the properties of the q-gamma function.

\subsection{The summation formula}\label{diofanto}
{}From \Ref{eq-SumFormula} we get 
\begin{equation}\label{sumform}
\Gamma(z,\tau,\sig)=
\exp{\left(-\frac i2\sum_{j=1}^\infty
\frac{\sin(\pi j(2z-\tau-\sig))}
{j\sin(\pi j\tau)\sin(\pi j\sig)}\right)}.
\end{equation}
The region of absolute convergence of this series also
include points where $\tau$ or $\sig$ do not have positive
imaginary part. If $\tau,\sig\in\C-\R$, then the series
converges absolutely if and only if
\begin{equation}\label{e-9022}
| \mathrm{Im}\,(2z-\tau-\sig)|
< | \mathrm{Im}\,(\tau)|+| \mathrm{Im}\,(\sig)|.
\end{equation}
Whenever both sides of the equations are in 
the region of convergence of this series we then clearly have
\[
\Gamma(z,-\tau,\sig)=\Gamma(\sig-z,\tau,\sig),\qquad
\Gamma(z,\tau,-\sig)=\Gamma(\tau-z,\tau,\sig).
\]
These formulae can be used to extend the definition
of  the elliptic gamma
function for $\sig,\tau\in\C-\R$, as we do in the next subsection.

\subsection{Extending the range of parameters}\label{s-33}
Since many operations we perform  do not preserve the upper half plane,
it is important to extend the range of values $\tau$ and $\sig$
can take. 
We set 
\begin{equation}\label{e-000}
(x;q^{-1})=\frac1{(qx;q)},\qquad (x;q^{-1},r)=\frac1{(qx;q,r)}
\qquad 
(x;q,r^{-1})=\frac1{(rx;q,r)}
\end{equation}
These formulae define a unique extension of the functions $(x;q)$,
$(x;q,r)$ to meromorphic functions on $\{(x,q,r)| |q|\neq1\neq |r|\}$
obeying \Ref{e-000}.
It is clear that the functional relations
\[
(qx,q)=\frac1{1-x}(x,q),\qquad (qx;q,r)=\frac1{(x;r)}(x;q,r),
\]
still hold in this larger domain.
Correspondingly, we extend the definition of $\theta_0$ and the
elliptic gamma function by using the same formulae in terms of the
infinite products. We obtain:
\[\theta_0(z,-\tau)=\frac1{\theta_0(z+\tau,\tau)},\quad
\Gamma(z,-\tau,\sig)=\frac1{\Gamma(z+\tau,\tau,\sig)},\quad
\Gamma(z,\tau,-\sig)=\frac1{\Gamma(z+\sig,\tau,\sig)}.
\]
A straightforward check gives the following result:
\begin{thm}\label{t-3}
The identities \Ref{EGF1} for $\theta_0$ 
hold for all $\tau\in\C-\R$. The
identities of Theorem \ref{t-1}  and Prop.\ \ref{p-3}
for $\Gamma$ and $\theta_0$ 
hold for all $z\in\C$,  $\tau,\sig\in\C-\R$ whenever both
sides are defined. The summation formula \Ref{sumform} is
valid for all $z\in\C$, $\tau,\sig\in\C-\R$ such that the sum converges
absolutely.
\end{thm}

However, the statements about the position of zeros and poles are
no longer valid.

\subsection{Approaching the real axis}
Here we notice that the series \Ref{sumform} actually
also converges for certain real values of $\tau$ or $\sig$.
Indeed let, for any $\alpha>1$, $X_\alpha$ denote the set of
real numbers $\tau$ such that 
$\mathrm{min}_{k\in\Z}|j\tau-k|>j^{-\alpha}$ for all
but finitely many integers $j>0$. By Khintchin's theorem (see, e.g., \cite{C},
Chapter VII), $X_\alpha$, and therefore also $X=\cup_{\alpha>1}  X_\alpha$,
 is the complement in $\R$ of a set of Lebesgue measure zero. The
set $X$ contains
in particular all irrational algebraic numbers.
For $\tau\in X$, one has, for some $\alpha>1$, the bound
$|\sin(\pi j\tau)|\geq j^{-\alpha}$ 
for all sufficiently large $j$.
Therefore, the series is also absolutely convergent if
\[
\tau\in X, \quad|\mathrm{Im}(2z-\sig)|<|\mathrm{Im}\,{}\sig|,
\]
or, by symmetry, if
\[
\sig\in X, \quad |\mathrm{Im}(2z-\tau)|<|\mathrm{Im}\,{}\tau|.
\]
These results may be summarized as follows.

 \medskip

\begin{proposition}\label{p-gamm} Let
$X=\cup_{\alpha>1}\{\tau\in\R\,|\,
\mathrm{min}_{k\in\Z}|j\tau-k|>j^{-\alpha},\,\forall j>>1\}$. Then
$N=\R-X$ has Lebesgue measure zero and
the series appearing in \Ref{sumform} converges absolutely for
all $\tau,\sig\in\C-N$ and $z\in\C$ obeying
\[
|\mathrm{Im}\,(2z-\tau-\sig)|
< |\mathrm{Im}\,(\tau)|+|\mathrm{Im}\,(\sig)|.
\]
\end{proposition}

A more precise result is the following ``wall crossing theorem'', 
which shows that the
values of $\Gamma(z,\tau,\sig)$ for real $\tau$ 
are obtained as suitable limits from either side.

\begin{thm}
\label{t-real} Let $X$ be the subset of the real line
of Prop.\ \ref{p-gamm} and suppose that  $\tau\in X$ and
 $|\mathrm{Im}(2z-\sig)|
<|\mathrm{Im}\,\sig|$. Then, as a function of $\epsilon\in\R$,
$\Gamma(z,\tau+i\epsilon,\sig)$, as given by the
convergent series \Ref{sumform}  is continuous at $\epsilon=0$.
\end{thm}

\begin{prf} We need to estimate  the terms in the
sum \Ref{sumform} uniformly in $\epsilon$. Let us assume
for definiteness that $\epsilon\geq 0$. The case $\epsilon\leq0$
is treated in the same way.  If $|\mathrm{Im}\,x|\geq \delta>0$, we have
$
c(\delta)\exp(|\mathrm{Im}\,x|)\leq|\sin\,x|\leq \exp(|\mathrm{Im}\,x|),
$
with $0<c(\delta)=1-e^{-2\delta}<1$
This implies the bound on the $j$th term of \Ref{sumform}:
\[
\left|\frac{\sin(\pi j(2z-\tau-i\epsilon-\sig))}
{j\sin(\pi j(\tau+i\epsilon))\sin(\pi j\sig)}
\right|
\leq c(\pi|\mathrm{Im}\,\sig|)^{-1}
\frac {e^{\pi j(|\mathrm{Im}(2z-\sig)-\epsilon|-|\mathrm{Im}\,\sig|)}}
{j\sin\,\pi j(\tau+i\epsilon)}.\]
Let $\alpha,N$ be such that 
$\mathrm{min}_{k\in\Z}|j\tau-k|>j^{-\alpha},\,\forall j\geq N$.
The next step is to find a lower bound for $\sin\,\pi j(\tau+i\epsilon)$
for $j\geq N$. This is done in two different ways depending on
whether $j\epsilon$ is small or large.

\noindent(a) If $\sinh\,\pi j\epsilon\,\leq j^{-\alpha}e^{\pi\epsilon j}$, the
triangle inequality can be used in the form
\begin{eqnarray*}
|\sin\,\pi j(\tau+i\epsilon)|
&=&
\frac12|
e^{\pi i\tau j-\pi\epsilon j}
-
e^{-\pi i\tau j+\pi\epsilon j}
|
\\
&\geq&
\frac12|
e^{\pi i\tau j}
-
e^{-\pi i\tau j}
|
e^{\pi \epsilon j}
-\frac12|
e^{\pi\epsilon j}
-
e^{-\pi\epsilon j}
|\\
&=&
|\sin\,\pi\tau j|e^{\pi\epsilon j}-\sinh\,\pi\epsilon j.
\end{eqnarray*}
Let $k\in\Z$ so that $|\tau j-k|\leq 1/2$. By using $|\sin\,\pi x|\geq
2|x|$ for ${-1/2\leq x\leq 1/2}$, we get $|\sin\,\pi \tau j|=
|\sin\,\pi(\tau j-k)|\geq 2j^{-\alpha}$. Thus
\[
|\sin\,\pi j(\tau +i\epsilon)|\geq 2
j^{-\alpha}e^{\epsilon \pi j}
-\sinh\, \pi\epsilon j\geq 
j^{-\alpha}e^{\epsilon \pi j}.
\]
(b) 
 If $\sinh\,\pi j\epsilon\,\geq j^{-\alpha}e^{\pi\epsilon j}$, the
triangle inequality implies
\begin{eqnarray*}
|\sin\,\pi j(\tau+i\epsilon)|
&=&
\frac12|
e^{\pi i\tau j-\pi\epsilon j}
-
e^{-\pi i\tau j+\pi\epsilon j}
|
\\
&\geq&
\frac12(
|e^{-\pi i\tau j+\pi\epsilon j}|
-
|e^{\pi i\tau j-\pi\epsilon j}|
)
\\
&=&
\sinh\,\pi\epsilon j\geq j^{-\alpha}e^{\pi\epsilon j}
\end{eqnarray*}
In both cases (a) and (b) we get the lower bound
\[
|\sin\,\pi j(\tau+i\epsilon)|\geq j^{-\alpha}e^{\pi\epsilon j},
\]
for all $j\geq N$.
Therefore we have the uniform bound on the $j$th term of the series
\[
\left|\frac{\sin(\pi j(2z-\tau-i\epsilon-\sig))}
{j\sin(\pi j(\tau+i\epsilon))\sin(\pi j\sig)}\right|
\leq c(\pi\mathrm{Im}\,\sig)^{-1}
j^{\alpha-1}
{e^{\pi j(|\mathrm{Im}(2z-\sig)|-|\mathrm{Im}\,\sig|)}}
\]
for all $j\geq N$. The sum over $j$ of this expression converges
if $|\mathrm{Im}(2z-\sig)|\leq|\mathrm{Im}\,\sig|$. So our
series is bounded, for  all $\epsilon$, by a
single absolutely convergent series.
It follows that the sum is a continuous function of $\epsilon$. 
\end{prf}

\section{Modular properties}\label{S-4}
We consider the transformation properties of the elliptic
gamma function under modular transformations of $\sig$
and $\tau$. We have the identities
\begin{thm}\label{t-4}
Suppose that $\tau,\sig,\sig/\tau,\tau+\sig\in\C-\R$. 
Let \begin{eqnarray*}
Q(z;\tau,\sig)&=&
\frac{z^3}{3\tau\sig}
-
\frac{\tau+\sig-1}{2\tau\sig}
z^2
+
\frac{\tau^2+\sig^2+3\tau\sig-3\tau-3\sig+1}
{6\tau\sig}
z
\\ & &
+
\frac1{12}
(\tau+\sig-1)
(\tau^{-1}+\sig^{-1}-1).
\end{eqnarray*}
Then
\begin{eqnarray}
\Gamma(z,\tau+1,\sig)&=&\Gamma(z,\tau,\sig+1)=
\Gamma(z,\tau,\sig),
\\
\Gamma(z,\tau+\sig,\sig)
&=&
\frac{\Gamma(z,\tau,\sig)}
{\Gamma(z+\tau,\tau,\sig+\tau)},
\\
\Gamma(
 z/\sig,\tau/\sig,-1/\sig)&=&
e^{i\pi Q(z;\tau,\sig)}
{\Gamma(
{({z-\sig})/\tau,-1/\tau,-\sig/\tau})}
{\Gamma(z,\tau,\sig)}
,\\
\Gamma(
 z/\tau,-1/\tau,\sig/\tau)&=&
e^{i\pi Q(z;\tau,\sig)}\label{e-u}
{\Gamma(
{({z-\tau})/\sig,-\tau/\sig,-1/\sig})}
{\Gamma(z,\tau,\sig)}
.
\end{eqnarray}
\end{thm}

\begin{prf}
We give the proof of  these identities in the domain
where the second and third arguments of  all gamma functions have 
 positive imaginary part, so that the gamma functions
are defined by the product formula. The general case
is reduced to this case by inserting the definitions of \ref{s-33},
as a straightforward check shows.
The first two identities are obvious, the third follows from the identity
\[(x;qr,r)(qx;q,qr)=(x;q,r),\] 
which is easy to check. The last
identity is obtained from \Ref{e-u} by exchanging
$\tau$ and $\sig$ and using the symmetry \Ref{eq-Gamma9}.

To prove \Ref{e-u}, we show that the ratio between the two sides of the
equation is a triply periodic meromorphic function and is therefore
constant, and determine the constant by evaluating the ratio at
a special value.
 
Let $A(z;\tau,\sig)$ be the ratio
\[
\frac{
\Gamma\left(\frac z\sig,\frac\tau\sig,-\frac1\sig\right)
}
{
\Gamma\left(\frac {z-\sig}\tau,-\frac1\tau,-\frac\sig\tau\right)
\Gamma(z,\tau,\sig)
}
\]
We have
\begin{eqnarray*}
\frac{A(z-1;\tau,\sig)}{A(z;\tau,\sig)}
&=&
\frac{
\theta_0\left(\frac z\sig,\frac\tau\sig\right)
}
{
\theta_0\left(\frac {z-\sig}\tau,-\frac\sig\tau\right)
}
\\
&=&
i
e^{\pi i(-\tau/{6\sig}+z/\sig-\sig/{6\tau}-z/\tau
-z^2/\tau\sig)}
\frac
{\theta_0\left(\frac {z}\tau,-\frac\sig\tau\right)}
{\theta_0\left(\frac {z-\sig}\tau,-\frac\sig\tau\right)}
\\
&=&
\exp(\pi iP(z;\tau,\sig)),
\end{eqnarray*}
with
\[
P(z,\tau,\sig)
=
-\frac{z^2}{\tau\sig}
+z(\tau^{-1}+\sig^{-1})-
\frac{\sig}{6\tau}-\frac\tau{6\sig}
-\frac12
=P(z,\sig,\tau).
\]
Similarly, we find
\begin{eqnarray*}
\frac
{A(z+\tau;\tau,\sig)}
{A(z;\tau,\sig)}
&=&
\exp(\pi iP(z;-1,\sig)),
\\
\frac
{A(z+\sig;\tau,\sig)}
{A(z;\tau,\sig)}
&=&
\exp(\pi iP(z;\tau,-1)).
\end{eqnarray*}
The polynomial $Q$ is designed to compensate for
these terms. This is most easily seen by setting
$Q(z;\tau,\sig,\rho)=Q(-z/\rho;-\tau/\rho,-\sig/\rho)$
which is symmetric in $\tau,\sig,\rho$ and
obeys
\[
Q(z+\tau;\tau,\sig,\rho)=
Q(z;\tau,\sig,\rho)+P(z;\sig,\rho).
\]
By 
using this identity for permutations of $\tau,\sig,\rho$, 
with $\rho=-1$, we deduce that 
$Ae^{-\pi i Q}$ is
a triply periodic meromorphic function of $z$ and
is thus constant. To compute this constant, we
set $z=(\tau+\sig-1)/2$. Then all gamma functions
are equal to one, and $Q(z;\tau,\sig)=0$. Thus
the constant is one.
\end{prf}

\section{Special values}\label{S-5}
Here we consider the degeneration of our three term
relations when the periods $\tau,\sig$ 
are linearly dependent. The simplest case is when
$\tau=\sig$. Then we can write the gamma function as
\[
\Gamma(z,\tau,\tau)=\prod_{j=0}^\infty
\left(
   \frac{1-e^{-2\pi iz}q^{j+2}}
        {1-e^{2\pi iz}q^j}
\right)^{j+1},\qquad q=e^{2\pi i\tau}.
\]
To express the result we need to recall
a simple property of the dilogarithm function.

\begin{proposition}\label{p-478}
 Let $\mathrm{Li}_2(x)=\sum_{j=1}^\infty
\frac {x^j}{j^2}$ be the dilogarithm and let for
$\mathrm{Im}\,t<0$,
\[
\psi(t)=\exp\left(t\ln(1-e^{-2\pi it})-\frac1{2\pi i}
\mathrm{Li}_2(e^{-2\pi it})\right),
\]
where the branch of the logarithm is determined
by $\ln(1-x)=-\sum_1^\infty x^j/j$, ($|x|<1$).
Then $\psi(t)$ has an analytic continuation
to a meromorphic function on the complex plane.
It has a zero of order $n$ at $t=n$ and a pole of order $n$ 
at $t=-n$ ($n=1,2,\dots$) and no other zeros or poles.
Moreover $\psi$ obeys the functional equation
\[
\psi(t+1)=(1-e^{-2\pi it})\psi(t),
\]
and the estimate
\[
\psi(t)=1+0(|\mathrm {Im}\,t|e^{-2\pi |\mathrm{Im}\,t|}),
\]
as $\mathrm{Im}\,t\to -\infty$.
\end{proposition}
\noindent{\it Proof:}
It is clear that the Taylor series defines a holomorphic function
on the lower half plane obeying the functional equation.
The singularities on the real axis can be studied using the
integral representation of the dilogarithm:
\[
\psi(t)=\exp\left(2\pi i\int_{t}^{-i\infty}\frac {s\,ds}{e^{2\pi is}-1}\right).
\]
This well-known formula may be checked by expanding the geometric series
in the integrand and integrating term by term.
{}From this formula we see that the only potential singularities of the
argument of the exponential function are at integer values
of $t$. At
$t=0$, however, the function is regular ($\psi(0)=\exp{i\pi/12}$).
The functional equation then implies the statement about zeros
and poles. In particular $\psi$ is single-valued at the
singularities, so that the integral representation defines
a meromorphic function with no other zeros or poles.

The estimate follows from the inequalities $|\ln(1-x)|\leq 2|x|$,
if $|x|$ is sufficiently small and
$|\mathrm{Li}_2(x)|\leq |x|\sum_{1}^\infty1/j^2$ if $|x|\leq1$.
$\square$

\begin{thm}\label{t-5} Let $\mathrm{Im}\,\tau>0$ and $z\in\C-(\Z+\tau\Z)$.
Then
\[
\Gamma(z,\tau,\tau)=
\frac
{e^{-\pi iQ(z;\tau,\tau)}}
{\theta_0\left(
\frac z\tau,-\frac1\tau
\right)}
\prod_{k=0}^{\infty}
\frac
{\psi\left(
\frac{k+1+z}\tau\right)}
{\psi\left(\frac{k-z}\tau\right)}\,.
\]
\end{thm}

The infinite product is convergent thanks to the estimate
of Prop.\ \ref{p-478}

The following calculation is not a completely rigorous
proof of this theorem, but it is more transparent than
the correct proof, which consists of showing that the
ratio between left and right-hand side is an entire doubly
periodic meromorphic function taking the value 1 at
a special point.

We start from the three term relation for $\Gamma$.
\[
\Gamma(z,\tau,\sig)=
e^{-\pi iQ(z;\tau,\sig)}
\frac
{\Gamma\left(\frac z\tau,-\frac1\tau,\frac \sig\tau\right)}
{\Gamma\left(\frac{z-\tau}\sig,-\frac1\sig,-\frac\tau \sig\right)}
=
\frac{e^{-\pi iQ(z;\tau,\sig)}}
{\theta_0\left(\frac z\sig,
-\frac{1}\sig\right)}
\frac{\Gamma\left(\frac z\tau,-\frac1\tau,\frac \sig\tau\right)}
{\Gamma\left(\frac{z}\sig,-\frac1\sig,-\frac\tau \sig\right)}\,.
\]
Let us take the limit $\sig\to\tau$. The limit of the
ratio of gamma function is delicate. Set $\sig=\tau(1+\epsilon)$ and
introduce multiplicative variables:
\[
q_1=e^{-\frac{2\pi i}\tau},
\quad
q_2=e^{-\frac{2\pi i}\sig},
\quad
r_1=e^{\frac{2\pi i\sig}\tau},
\quad
r_2=e^{-\frac{2\pi i\tau}\sig},
\quad
x_1=e^{\frac{2\pi iz}\tau},
\quad
x_2=e^{\frac{2\pi iz}\sig}.
\]
Then, by the summation formula,
\[
\ln\frac
{\Gamma\left(\frac z\tau,-\frac1\tau,\frac \sig\tau\right)}
{\Gamma\left(\frac{z}\sig,-\frac1\sig,-\frac\tau \sig\right)}
=
-
\sum_{j=1}^\infty
\frac
{
(q_1r_1/x_1)^j-x_1^j
}
{
j(1-q_1^j)(1-r_1^j)
}
+
\sum_{j=1}^\infty
\frac
{
(q_2r_2/x_2)^j-x_2^j
}
{
j(1-q_2^j)(1-r_2^j)
}
\,.
\]
We now expand the various terms around $\epsilon=0$. We
get
\[
q_2=q_1(1+\frac{2\pi i}{\tau}\epsilon)+O(\epsilon^2),
\qquad
r_1=1+2\pi i\epsilon+O(\epsilon^2),
\]
\[
r_2=1+2\pi i\epsilon+O(\epsilon^2),
\qquad
x_2=e^{\frac{2\pi iz}\tau}(1-2\pi i\frac z\tau\epsilon).
\]
The singular terms have the expansion
\[
\frac{1}{1-r_1^j}=-\frac1{2\pi ij\epsilon}+\frac12+O(\epsilon)\,,
\]
\[
\frac{1}{1-r_2^j}=-\frac1{2\pi ij\epsilon}-\frac1{2\pi ij}+
\frac12+
O(\epsilon)\,.
\]
Inserting this in the summation formula yields
\begin{eqnarray*}
\ln\frac
{\Gamma\left(\frac z\tau,-\frac1\tau,\frac \sig\tau\right)}
{\Gamma\left(\frac{z}\sig,-\frac1\sig,-\frac\tau \sig\right)}
=-
\sum_1^\infty
\frac
{q_1^jx_1^{-j}\left(\frac1\tau+\frac z\tau\right)}
{j(1-q_1^j)}
-
\sum_1^\infty
\frac
{x_1^{j}}
{j(1-q_1^j)}
\cdot
\frac z\tau\\
-\sum_1^\infty
\frac
{q_1^jx_1^{-j}-x_1^j}
{2\pi i\,j^2(1-q_1^j)}
-
\sum_1^\infty
\frac
{q_1^jx_1^{-j}-x_1^j}
{j(1-q_1^j)^2}
\cdot\frac {q_1^j}\tau+O(\epsilon).
\end{eqnarray*}
After expanding the denominators  into geometric
series and exchanging the summations, the sums over $j$ become
Taylor series for (di)logarithms. The result is the formula
of the theorem.

\begin{corollary}\label{c-4} Let $\tau\to 0$ on a ray
$\{s\tau_0\,|\,s>0\}$ with $\mathrm{Im}(\tau_0)>0$ and let $z=u+v\tau_0$
with $-1<u<0$, $v\in\R$ be fixed. Then
\[
\Gamma(z,\tau,\tau)=e^{-\pi iQ(z;\tau,\tau)}(1+\mathcal{O}(e^{-\frac c
{\mathrm{Im}\,\tau}}))
\]
for some $c>0$ depending on $z$, $\tau_0$.
\end{corollary}

\noindent{\bf Remark.} As $\tau\to 0$ along this ray, the zeros and poles
of $\Gamma(z,\tau,\tau)$ as a function of $z$ accumulate on
the lines $n+s\tau_0$, $n\in\Z$, $s\in\R$. The assumption on 
$z$ means that $z$ lies between two lines. One can relate this
case to the more general case of $z$ between any two other lines
by using the fact that $\Gamma(z,\tau,\tau)$ is 1-periodic.
\medskip

\noindent{\it Proof of Corollary \ref{c-4}:}
The assumption on $z$ implies that the arguments of the
$\psi$ functions in Theorem \ref{t-5}  obey
\[
\mathrm{Im}\frac{k+1+z}{\tau}=
({k+1}+u)\mathrm{Im}\frac1{\tau},\qquad
\mathrm{Im}\frac{k-z}{\tau}=
(k-u)\mathrm{Im}\frac1{\tau},
\]
Since $\mathrm{Im}(1/\tau)\to-\infty$ and $k+1+u, k-u>0$ for
all $k=0,1,2,\dots$, we can use the estimate
of Prop.\ \ref{p-478} to show that the product of ratios of
$\psi$ functions tends to $1$ as $\tau\to0$ with an error term that is 
smaller than $c_1e^{-c/\mathrm{Im}\,\tau}$ for some constants
$c_1,c>0$. Similarly $\theta_0$ tends to one, and we are left
with the exponential of $Q$. $\square$

\medskip
 More generally, one can find similar formulae when the
periods $\tau,\sig$ are linearly dependent over
the rationals. The following two theorems reduce (in different ways)
this computation
to the case studied above.

\begin{thm}\label{t-6}
Let $a,b$ be positive integers. Then
$$
\Gamma(z,a\tau,b\tau)
=\prod_{r=0}^{b-1}\prod_{s=0}^{a-1}
\Gamma(z+(ar+bs)\tau,ab\tau,ab\tau).
$$
\end{thm}
\begin{prf}
Let $q=e^{2\pi i\tau}$, $x=e^{2\pi i z}$. Then
\[
\Gamma(z,a\tau,b\tau)=
\frac {(q^{a+b}x^{-1};q^a,q^b)}{(x;q^a,q^b)}\,.
\]
We first prove
an identity for double products
\begin{eqnarray*}
(x;q^a,q^b)&=&
\prod_{j,k=0}^\infty(1-xq^{aj+bk})\\
&=&
\prod_{r=0}^{b-1}\prod_{s=0}^{a-1}
\prod_{j,k=0}^\infty(1-xq^{a(r+bj)+b(s+ak)})\\
&=&
\prod_{r=0}^{b-1}\prod_{s=0}^{a-1}
(xq^{ar+bs};q^{ab},q^{ab}).
\end{eqnarray*}
If we replace $x$ by $q^{a+b}x^{-1}$ in this identity and change variables
$r\to b-1-r$, $s\to a-1-s$, we obtain
\[
(q^{a+b}x^{-1};q^a,q^b)
=
\prod_{r=0}^{b-1}\prod_{s=0}^{a-1}
(x^{-1}q^{2ab-ar-bs};q^{ab},q^{ab}).
\]
Taking the ratio we get the desired identity for gamma functions.
\end{prf}

\medskip

\noindent{\bf Examples.}
\begin{enumerate}
\item
$
\Gamma(z,\tau,3\tau)
=
\Gamma(z,3\tau,3\tau)
\Gamma(z+\tau,3\tau,3\tau)
\Gamma(z+2\tau,3\tau,3\tau).
$
\item
$
\Gamma(z,2\tau,3\tau)
=
\prod_{j\in\{0,2,3,4,5,7\}}\Gamma(z+j\tau,6\tau,6\tau)
$
\item
Setting $a=b$ and rescaling $\tau$ we get
\[
\Gamma(z,\tau,\tau)=\prod_{j=0}^{2(a-1)}\Gamma(z+j\tau,a\tau,a\tau)^{
a-|j-a+1|}.
\]
\end{enumerate}
Inserting the formula of Theorem \ref{t-5} into the formula
of Theorem \ref{t-6} yields an expression for $\Gamma(z,a\tau,b\tau)$
in terms of dilogarithms.
As an application we can compute from Corollary \ref{c-4},
the asymptotics of the infinite products
\[
\Gamma(z,a\tau,b\tau)=\prod_{j=0}^\infty
\left(
\frac{1-q^{j+a+b}e^{-2\pi iz}}{1-q^je^{2\pi iz}}\right)^{N_{a,b}(j)}.
\]
Here $N_{a,b}(j)$ denotes the number of ways $j$ can be
written as $j=ar+bs$ with $r,s$ nonnegative integers.
We have, as $\tau\to 0$ as in Corollary \ref{c-4},
\[
\Gamma(z,a\tau,b\tau)=e^{-\pi i\sum_{r=0}^{b-1}
\sum_{s=0}^{a-1}
Q(z+(ar+bs)\tau;ab\tau,ab\tau)}(1+\mathcal{O}(e^{-\frac c
{\mathrm{Im}\,\tau}})).
\]

\begin{thm}\label{formula}
Let the greatest common divisor of natural numbers $a, b$ be 1.
Consider the function
\[
\Gamma(z,a\tau,b\tau)=\prod_{j=0}^\infty
\left(
\frac{1-q^{j+a+b}e^{-2\pi iz}}{1-q^je^{2\pi iz}}\right)^{N_{a,b}(j)}
\]
where $N_{a,b}(j)$ denotes the number of ways $j$ can be
written as $j=ar+bs$ with $r,s$ nonnegative integers.
Then
\[
\Gamma(z,a\tau,b\tau)^{ab}\,=\,\Gamma(z,\tau,\tau)
\,\prod_{k=0}^{ab-1}\theta_0(z+k\tau, ab\tau)^{\alpha_k}
\]
where $\alpha_k=-ab+k+1$ if $k=ar+bs$ for some integers $r,s\geq 0$
and
$\alpha_k=k+1$ if $k$ cannot be represented in this form.

\end{thm}

\noindent{\bf Example.}
\begin{eqnarray*}
&\Gamma(z,2\tau,3\tau)^6=\Gamma(z,\tau,\tau)\times
\notag
\\
&\theta_0(z,6\tau)^{-5}\theta_0(z+\tau,6\tau)^{2}
\theta_0(z+2\tau,6\tau)^{-3}
\theta_0(z+3\tau,6\tau)^{-2}
\theta_0(z+4\tau,6\tau)^{-1}\,.
\notag
\end{eqnarray*}

To prove the Theorem we shall use the following Lemma.

\begin{lemma} Let $k\in\{0,...,ab-1\}$.
\begin{enumerate}
\item[1.] If $k=ar+bs$ for some integers $r,s\geq 0$, then
$ab-k$ cannot have the form $ab-k=a(i+1)+b(j+1)$ for
some integers $i,j\geq 0$.
\item[2.] If $k$ cannot be represented in the form
$k=ar+bs$ for some integers $r,s\geq 0$, then
$ab-k=a(i+1)+b(j+1)$ for some integers $i,j\geq 0$.
\end{enumerate}
\end{lemma}

\noindent{\it Proof of the Lemma:} If $k=ar+bs$ and
$ab-k=a(i+1)+b(j+1)$, then $ab=a(r+i+1)+b(s+j+1)$.
Since $a,b$ are relatively prime, this leads to a contradiction.
Part 1 is proved.

If $k$ cannot be represented in the form
$k=ar+bs$ for some integers $r,s\geq 0$, then
$k$ can be represented in the form $k=ar'+bs'$
where $0<r'<b,\,{}$ $ s'<0$.
This gives the desired representation for $ab-k$,
 $ab-k= a(b-r')-bs'$. Part 2 is proved.

\noindent{\it Proof of Theorem \ref{formula}:}
We have
\[
\Gamma(z,a\tau,b\tau)^{ab}=
\prod_{k=0}^{ab-1}
\prod_{s=0}^\infty
\frac{(1-q^{ab(s+1)}q^{-k}e^{-2\pi iz})^{ab\beta_{k,s}}}
{(1-q^{abs} q^ke^{2\pi iz})^{ab\gamma_{k,s}}}
\]
where $\beta_{k,s}$ is the number of ways $ab(s+1)-k$
can be written as $ab(s+1)-k = a(i+1)+ b(j+1)$ with nonnegative integers $i,j$
and $\gamma_{k,s}$ is the number of ways $abs+k$
can be written as $abs+k = ai+ bj$ with nonnegative integers $i,j$.
It is easy to see that $\beta_{k,s}=s+\beta_{k,0}$
and $\gamma_{k,s}=s+\gamma_{k,0}$. By the Lemma, $\beta_{k,0}+\gamma_{k,0}=1$.
Notice also that $\alpha_k = k+1 - ab\gamma_{k,0}$.

Thus we have
\begin{eqnarray*}
&\Gamma(z,a\tau,b\tau)^{ab}=
\prod_{k=0}^{ab-1}
\prod_{s=0}^\infty
\frac{(1-q^{ab(s+1)}q^{-k}e^{-2\pi iz})^{ab(s+1) - ab\gamma_{k,0}}}
{(1-q^{abs} q^ke^{2\pi iz})^{abs + ab \gamma_{k,0}}}
=
\notag  
\\
&\prod_{k=0}^{ab-1}
\left( \prod_{s=0}^\infty
(1-q^{ab(s+1)}q^{-k}e^{-2\pi iz})(1-q^{abs} q^ke^{2\pi iz})\right) ^{k+1-ab\gamma_{k,0}}
\,{}\,\times
\notag
\\
&
\prod_{k=0}^{ab-1}
\prod_{s=0}^\infty
\frac{ (1-q^{ab(s+1)}q^{-k}e^{-2\pi iz})^{ab(s+1) -k-1} }
{(1-q^{abs} q^ke^{2\pi iz})^{abs + k+1}}=
\notag
\\
& \Gamma(z,\tau,\tau)
\,\prod_{k=0}^{ab-1}\theta_0(z+k\tau,ab\tau)^{\alpha_k}.
\end{eqnarray*}
The Theorem is proved.

\section{The phase function and the semiclassical limit}\label{S-6}

Here we introduce the ``phase function'', which is the ratio
of elliptic gamma functions appearing in hypergeometric
integrals \cite{FTV}. It obeys identities which are direct
consequences of the identities for elliptic gamma functions.
We discuss the phase function here since it has a semiclassical
limit in which the identities reduce to more familiar differential
equations and modular properties of theta functions.

We keep the notation of \ref{s-egf}
and introduce a new variable $a$,  setting
$\alpha=e^{2\pi ia}$. The phase function is defined as the ratio
\begin{equation}\label{e-O}
\Omega_a(z,\tau,\sig)
=\frac{\Gamma(z+a,\tau,\sig)}
{\Gamma(z-a,\tau,\sig)}
=\frac{(qr/x\alpha;q,r)(x/\alpha;q,r)}
{(x\alpha;q,r)(qr\alpha/x;q,r)}.
\end{equation}
The following identities are direct translations of identities
for gamma functions:
\begin{eqnarray}
\Omega_a(z+\sig,\tau,\sig)
&=&\frac{\theta_0(z+a,\tau)}
{\theta_0(z-a,\tau)}
\Omega_a(z,\tau,\sig),\label{eo-1}\\
\Omega_a(z+\tau,\tau,\sig)
&=&\frac{\theta_0(z+a,\sig)}
{\theta_0(z-a,\sig)}
\Omega_a(z,\tau,\sig),\label{eo-2}\\
\Omega_a(z+1,\tau,\sig)
&=&
\Omega_a(z,\tau,\sig),\label{eo-3}
\end{eqnarray}
The modular properties of this function also follow directly from those
of the gamma function. In particular
\begin{eqnarray}
\Omega_a(z,\tau+\sig,\sig)
&=&
\frac{\Omega_a(z,\tau,\sig)}
{\Omega_a(z+\tau,\tau,\sig+\tau)},\label{eo-4}
\\
\Omega_a(z,\tau+1,\sig)
&=&
\Omega_a(z,\tau,\sig),\label{eo-5}\\
\Omega_{\frac a\tau}
(\textstyle{
{\frac z\tau},{\frac\sig\tau},{-\frac1\tau}})
&=&e^
{\pi i(Q(z+a;\tau,\sig)-Q(z-a;\tau,\sig)}
\Omega_{\frac a\sig}(\textstyle{
{\frac{z-\tau}\sig},{-\frac1\sig},{-\frac\tau\sig}
})
\Omega_a(z,\tau,\sig).\label{eo-6}
\end{eqnarray}
The argument of the  exponential function is the polynomial
\begin{equation}\label{e-99}
R_a(z,\tau,\sig)=\frac{\pi ia}{3\tau \sig}(6z^2-6(\tau+\sig-1)z+
2a^2+\tau^2+\sig^2+3\tau \sig-3\tau-3\sig+1).
\end{equation}
The summation formula \Ref{sumform} for $\Gamma$ implies
\begin{equation}\label{e-26}
\Omega_a(z,\tau,\sig)=\exp
\left(-i\sum_{\ell=1}^\infty
\frac{\cos(\pi \ell(2z-\tau-\sig))
\sin(2\pi \ell a)}{\ell\sin(\pi\ell\tau)\sin(\pi\ell \sig)}\right).
\end{equation}
By using this formula we can compute the {\em semiclassical limit}
$\sig\to 0$, $a\to 0$ with $\beta=2a/\sig$ fixed: if $0<\mathrm{Im}\,z<\mathrm{Im}\,\tau$
we are in the region of convergence of the series \Ref{e-26} and we get
\begin{eqnarray*}
u(z,\tau)
&=&
\lim_{\epsilon\to 0}
\Omega_{\epsilon}(z,\tau,2\epsilon/\beta)
\\
&=&
\exp\left(-i\sum_{\ell=1}^\infty\frac{\cos(\pi\ell(2z-\tau))}{\ell\sin(\pi\ell\tau)}
\beta\right)
\\&=&
\theta_0(z,\tau)^\beta, \qquad 0<\mathrm{Im}\,z<\mathrm{Im}\,\tau,
\end{eqnarray*}
cf.\ \Ref{e-tsf}.
To avoid discussing cuts we assume here that $\beta$ is an integer.

Now the identities for $\Omega$ become differential and difference
equations for the limit $u(z,\tau)$. Eq.\ \Ref{eo-1} is the obvious
differential equation $u'(z,\tau)=\beta(\theta'_0(z,\tau)/\theta_0(z,\tau))\,u(z,\tau)$,
where the derivative with respect to $z$ is denoted by a prime.
By \Ref{e-modt0} we have in the semiclassical limit
\[
\frac{\theta_0(z+\epsilon,\sig)}
{\theta_0(z-\epsilon,\sig)}\to e^{-\pi i\beta (2z-1)}.
\]
Then from \Ref{eo-2}, \Ref{eo-3} we see that the semiclassical limit exists
for almost all $z$ and has the theta function property
\[
u(z+\tau,\tau)=-e^{-\pi i\beta (2z-1)}u(z,\tau),\qquad
u(z+1,\tau)=u(z,\tau).
\]
Therefore the limit is $\theta_0^\beta$ for almost all $z,\tau$.
Let us now consider the semiclassical limit of \Ref{eo-4}.
Expanding \Ref{e-26} yields
 $\Omega_{\epsilon}(z+\tau,\tau,\tau+2\epsilon/\beta)
 =  1-2\pi i\epsilon \,r(z,\tau)+O(\epsilon^2)$, with
\begin{eqnarray*}
r(z,\tau)=\sum_{j=1}^\infty
  \frac{\cos(2\pi jz)}{\sin^2(\pi j\tau)}\\
\end{eqnarray*}
Therefore both sides of \Ref{eo-5} tend to the same limit. But
the terms of order $\epsilon$ reduce to the differential equation
\[ \frac{\partial u}{\partial\tau}=\pi i\beta\, r(z,\tau)\,u.
\]
The identity \Ref{eo-5} becomes $u(z,\tau+1)=u(z,\tau)$.
We now turn to \Ref{eo-6}.
Let us assume that $\mathrm{Im}\,\sig$, $\mathrm{Im}\,\sig/\tau>0$,
$\mathrm{Im}\,\tau>0$. Then all factors  in the infinite products
\Ref{e-O} tend to one in the semiclassical limit and we get
\[
\Omega_{\frac a\sig}(\textstyle{
{\frac{z-\tau}\sig},{-\frac1\sig},{-\frac\tau\sig}
})\to 1.
\]
Therefore,  \Ref{eo-6} implies in the limit the modular
transformation properties of $u$:
\[
u(z/\tau,-1/\tau)=e^{R_0(z,\tau)}u(z,\tau).
\]
The expression $R_0$ is (see \Ref{e-99})
\[
R_0(z,\tau)=\lim_{\epsilon\to 0}R_\epsilon(z,\tau,2\epsilon\beta^{-1})
=\pi i\beta
\left(
 \frac{z^2}\tau
-z
+\frac{z}\tau
+\frac\tau6
-\frac12
+\frac1{6\tau}
\right).
\]
in agreement with \Ref{e-modt0}.

\section{A cohomological interpretation}\label{s-last}
Here we give an interpretation of the modular identities
obeyed by elliptic gamma functions. These identities may
be formulated in terms of the cohomology of $\mathrm{SL}(3,\Z)\semidirect\Z^3$
and automorphic forms ``of degree 1''.

\subsection{Automorphic forms}
Let us first review the well-known theory in degree zero. 
First some notational preliminaries.
We write all groups multiplicatively unless stated otherwise.
We denote by $[u]$ the equivalence class of an element $u$
in a quotient of abelian groups.
If $G$ is a group, a $G$-module $A$ is an abelian group with
a group homomorphism $\rho:G\to\mathrm{Aut}(A)$. The group $C^j(G,A)$ 
of $j$-cochains is the group of maps $\phi:G^j\to A$,
such that $\phi(g_1,\dots,g_j)=1$ if some $g_i=1$. One sets
$C^0(G,A)=A$. The
differential $\delta=\delta_j:C^j(G,A)\to C^{j+1}(G,A)$ is defined by
\begin{eqnarray*}
\lefteqn{\delta\phi(g_1,\dots,g_{j+1})=[\rho(g_1)\phi(g_2,\dots,g_{j+1})}\\
&&\prod_{i=1}^j\phi(g_1,\dots,
g_ig_{i+1},\dots,g_{j+1})^{(-1)^i}]^{(-1)^{j+1}}\phi(g_1,\dots,g_j),
\end{eqnarray*}
for $j\geq1$ and $\delta_0\phi(g)=\phi/\rho(g)\phi$ for $j=0$.
The $j$-th cohomology group of $G$ with coefficients $A$ is then
$H^j(G,A)=\mathrm{Ker}\,\delta_j/\mathrm{Im}\,\delta_{j-1}$  for 
$j\geq 1$ and $\mathrm{Ker}\,\delta_0$ for $j=0$.

Suppose  now
that $X$ is a connected complex manifold with nice action
of a group $G$, so that $X/G$ is a complex manifold.
Let $N$ be the multiplicative group  of the field of meromorphic
functions on $X$ and $M$ be the subgroup of nowhere vanishing
holomorphic functions. These groups are $G$-modules,
i.e., we have homomorphisms $G\to \mathrm{Aut}(N)$,
$G\to \mathrm{Aut}(M)$: $g\in G$ is mapped to the
automorphism $u\mapsto u(g^{-1}\cdot)$. The group of invariants
$H^0(G,N)=N^G$ is identified with the group of non-zero meromorphic
functions on $X/G$. 
Now let $\phi:G\to M$ be a 1-cocycle with coefficients $M$. 

We define\footnote{We consider meromorphic invertible automorphic forms. In particular
$0$ is not considered as an automorphic form and the additive structure is disregarded}
 {\em automorphic forms of type $\phi$} to be
functions $u\in N$ so that 
\[
u(x)=\phi(g,x)u(g^{-1}x),
\]
or 
$\phi=\delta u$ in $C^1(G,N)$. The 1-cocycle $\phi$ is
called the factor of automorphy of the automorphic form $u$. 

Automorphic forms corresponding to cocycles $\phi_1$, $\phi_2=\phi_1\delta\psi$, ($\psi\in C^0(G,M)$)
in the same cohomology class in $H^1(G,M)$ are in one-to-one correspondence via
$u\mapsto \psi u$. Thus it is convenient to consider equivalence classes of
 automorphic forms modulo $M$, which are associated to a cohomology
class of 1-cocycles: for each $[\phi]\in H^1(G,M)$, an {\it automorphic
class} of type $[\phi]$ is a class $[u]\in (N/M)^G=H^0(G,N/M)$ so that 
$[\delta u]=[\phi]$.

The basic properties of automorphic classes can be
expressed as follows: to any
short exact sequence $1\to M\to N\to N/M\to 1$ of $G$-modules
is associated the long exact sequence of cohomology
groups
\[
\cdots\to H^j(G,M)
\stackrel{i_*}\rightarrow
 H^j(G,N)
\stackrel{p_*}\rightarrow
 H^j(G,N/M)
\stackrel{\delta_*}{\rightarrow} H^{j+1}(G,M)
\stackrel{i_*}\rightarrow\cdots,
\]
The image of a class $[u]$ by the connecting homomorphism $\delta_*$ is obtained 
by viewing any representative $u$ as a cochain $u\in C^j(G,N)$, 
and setting $\delta_*[u]=[\delta u]$,

In our case, with $j=0$, the set of automorphic classes
of type $[\phi]\in H^1(G,M)$ is $\delta_*^{-1}[\phi]\subset H^0(G,N/M)$.
Exactness at $H^1(G,M)$ tells
us that the factors of automorphy $[\phi]$ for which there exist
automorphic forms are those in the kernel of $i_*$ and exactness at $H^0(G,N/M)$
implies
that the group $H^0(G,N)$ of non-zero meromorphic functions on $X$ acts
transitively on  $\delta_*^{-1}[\phi]$.

For example, take $G$ to be the free abelian group on two
generators $t_1,t_2$ acting on $X=\C$ by $t_1z=z+1$, 
$t_2z=z+\tau$. Then $H^0(G,N)$ is the group of non-zero
elliptic functions, $\theta_0$ represents a class of
$H^0(G,N/M)$ and $[\phi]=\delta_*[\theta_0]$ is the class
of the 1-cocycle
\[
\phi(t_1^lt_2^m,z)=\frac{\theta_0(z)}{\theta_0(z-l-m\tau)}
                  =e^{-\pi im(2z+1-(m+1)\tau)}.
\]
\newcommand{\abcd}{\textstyle{\left(\begin{array}{cc}a&b\\ c&d\end{array}\right)}}
More generally, let $H$ be the upper half plane
and let $G=\mathrm{SL}(2,\Z)\semidirect\Z^2$,
 the ``Jacobi group'',
act on $X=\C\times H$ by $(\abcd,\vec n)(z,\tau)=
((z+n_1+n_2\tau)/(c\tau+d),(a\tau+b)/(c\tau+d))$.
Then one has 
\[
\theta_0(y)=
\phi(g,y)\theta_0(g^{-1}y),\qquad y\in \C\times H,\quad g\in G.
\]
The factor of automorphy is defined on generators $t_1,t_2$ of $\Z^2$
and
$S={\scriptstyle\left(\begin{array}{cc}0&1\\ -1&0\end{array}\right)}$,
$T=\left(\begin{array}{cc}1&1\\ 0&1\end{array}\right)$ of $\mathrm{SL}
(2,\Z)$ by
\[
\phi(t_2,z,\tau)=e^{-\pi i(2z-2\tau+1)},
\qquad
\phi(S,z,\tau)=e^{\pi i\left(
-\frac{z^2}{\tau}
+z
-\frac{z}{\tau}
-\frac{\tau}{6}
-\frac{1}{6\tau}
+\frac{1}{2} \right)},
\]
and $\phi(t_1,z,\tau)=\phi(T,z,\tau)=1$. The value of $\phi$
on arbitrary group elements is then uniquely determined by the
cocycle relation.

Thus $\theta_0$ is a representative of an automorphic class
of type $[\phi]\in H^1(G,M)$. 

Let us review the proof that $[\phi]$ (and thus
$[\theta_0]\in H^0(G,N/M)$) is a non-trivial
cohomology class  by showing that the corresponding
first Chern class is non-trivial. The first Chern class
$c_1([\phi])$ is the image of $[\phi]$ under the
connecting homomorphism $H^1(G,M)\to H^2(G,2\pi i\Z)$ 
in the long exact sequence associated to $0\to 2\pi i\Z
\to\mathcal O(\C\times H)\stackrel{\exp}{\to}
 M\to 1$. Here $\mathcal O(X)$ is the additive
group of holomorphic functions on $X$. Then one checks
that $c_1([\phi])$ is non-trivial by showing that 
it is sent to a generator of $H^2(\Z^2,2\pi i\Z)\simeq \Z$
by the map $H^2(G,2\pi i\Z)\to H^2(\Z^2,2\pi i\Z)$
induced by inclusion. We do this calculation explicitly in
degree 1 below.

\subsection{Automorphic forms of degree 1}
There is an obvious generalization of these constructions
one degree higher: so for $G$-modules $M\subset N$ consisting
of certain classes of functions on an $G$-set $X$ we consider
the piece
\[
 H^1(G,N)\stackrel{p_*}\rightarrow H^1(G,N/M)\stackrel{\delta_*}\rightarrow H^2(G,M)
\stackrel{i_*}\rightarrow H^2(G,N)
\]
of the long exact sequence. A {\em degree 1 automorphic class} of type $[\phi]
\in \mathrm{Ker}(i_*:H^2(G,M)\to H^2(G,N))$ is then a class in $\delta_*^{-1}[\phi]$. 
Elements of these equivalence classes we call
{\em degree 1 automorphic forms}.
Degree 1 automorphic classes
 of type $[\phi]$ are acted upon transitively by $H^1(G,N)$.

We wish to show that $\theta_0,\Gamma$ are values on generators of a
degree 1 automorphic form for certain modules $M,N$ over $\mathrm{SL}(3,\Z)
\semidirect
\Z^3$.
As above let us start by considering the translation subgroup $\Z^3$.

\begin{proposition}\label{p-sic}
 Fix $\sigma,\tau$ in the upper half plane.
Let  $G$ be the free abelian group on three generators
$t_1,t_2,t_3$ acting on $\C$ by $t_1z=z+1$,
$t_2z=z+\sigma$, $t_3z=z+\tau$.
Let $N$ be the group of non-zero meromorphic functions
on $\C$, and $M$ the subgroup on nowhere vanishing holomorphic
functions. Then 
\[
u(t_1^lt_2^mt_3^n,z)=\prod_{j=1}^m\theta_0(z-j\sigma,\tau)
\]
represents a non-trivial class $[u]$ in $H^1(G,N/M)$. It corresponds to
the class $\delta_*[u]\in H^2(G,M)$ of the cocycle
\[
\phi(t_1^lt_2^mt_3^n,t_1^{l'}t_2^{m'}t_3^{n'},z)
=
e^{\pi i(nm'(2z+1)-m'n(n+1)\tau-nm'(m'+1+2m)\sigma)}.
\]
\end{proposition}
\medskip
\noindent{\it Proof:}
We need to compute 
\[
\delta u(t_1^lt_2^mt_3^n,t_1^{l'}t_2^{m'}t_3^{n'},z)
=\frac
{\prod_{j=1}^{m'}\theta_0(z-l-m\sigma-n\tau-j\sigma,\tau)
\prod_{j=1}^{m}\theta_0(z-j\sigma,\tau)}
{\prod_{j=1}^{m+m'}\theta_0(z-j\sigma,\tau)}
\,.\]
Using the transformation properties of $\theta_0$ under
translations by $\Z+\tau\Z$ we see that this coboundary
is in $M$, so that $u$ represents a class in $H^1(G,N/M)$.
The class of $\delta u$ in $H^2(G,M)$ is $\delta_*[u]$
and is easily computed to give the above expression.

To show that the class of $u$ is non-trivial it is
sufficient to show that $[\phi]$ is non-trivial.
To this purpose we compute the analogue of the first
Chern class: let $\mathcal O(\C)$ denote
the additive group of holomorphic functions on $\C$.
To the short exact sequence $0\to2\pi i\Z\to
\mathcal O(\C)\stackrel{\exp}\rightarrow M\to 1$ of
$G$-modules (with trivial action on $2\pi i\Z$) there
corresponds a long exact sequence and in particular
a connecting homomorphism
\[
c_1:H^2(G,M)\to H^3(G,2\pi i\Z)\simeq\Z
\]
We claim that $[\phi]$ is mapped to a generator under
this homomorphism and thus is non-trivial. The calculation
goes as follows: if we write $\phi(g_1,g_2,z)=\exp(R(g_1,g_2,z))$,
with $R$ the polynomial appearing in the exponential function
above,
then, according to the rules of homological algebra,
\[
c_1(\phi)(g_1,g_2,g_3)=
-R(g_2,g_3,g_1^{-1}z)+R(g_1g_2,g_3,z)
-R(g_1,g_2g_3,z)+R(g_1,g_2,z).
\]
If $g_j=t_1^{l_j}t_2^{m_j}t_3^{n_j}$, we then obtain
$c_1(\phi)(g_1,g_2,g_3)=2\pi il_1n_2m_3$ which is indeed
the class of a generator of $H^3(G,2\pi i\Z)$ (see \cite{MacL}
Section VI.6).
$\square$

Let us now give a similar interpretation for the
elliptic gamma function. By an {\em invertible analytic
function} on a complex manifold $X$ we mean an equivalence
class of pairs $(f,D)$ where $D\subset X$ is a dense open 
subset and $f$ is a holomorphic, nowhere vanishing function
on $D$. Two pairs $(f_1,D_1)$, $(f_2,D_2)$ are equivalent
if $f_1=f_2$ on $D_1\cap D_2$. Invertible analytic
functions form a group with respect to the  pointwise product
$(f_1,D_1)\cdot(f_2,D_2)=(f_1f_2,D_1\cap D_2)$.

It is convenient to pass from affine coordinates $\tau,\sigma$
to homogeneous coordinates $x_1,x_2,x_3$.
Let $G=\mathrm{SL}(3,\Z)\semidirect \Z^3$ act on $X=\C\times \C^3$
by $(A,\vec n)(z,\vec x)=(z+\vec n\cdot\vec x,A\vec x)$,
$A\in\mathrm{SL(3,\Z)}$, $\vec n=(n_1,n_2,n_3)\in\Z^3$,
$\vec x\in\C^3$, $z\in \C$.

The group $G$ has generators $e_{i,j}$, $1\leq i,j\leq3, i\neq j$
and $t_i$, $1\leq i\leq 3$. The elementary matrix $e_{i,j}$ is
the element of $\mathrm{SL}(3,\Z)$ which differ from the identity
matrix by having the $i,j$ matrix element equal to $1$. The $t_i$
are the canonical generators of the $\Z^3$ subgroup.

\begin{thm}\label{t-co}
 Let $G=\mathrm{SL}(3,\Z)\semidirect\Z^3$ act on $X=\C\times\C^3$
as above.
Let $N$ be the 
$G$-module of invertible analytic functions on $X$ such that
$f(\lambda z,\lambda \vec x)=f(z,\vec x)$ for all $\lambda\in\C-\{0\}$.
Let $M$ be the submodule of functions of the form
$\exp 2\pi i f$ with $f\in\Q(x_1/x_3,x_2/x_3)[z/x_3]$ a polynomial
in $z$ with coefficients in the rational functions of $\vec x$.
Then the classes in $N/M$ of the functions
\begin{eqnarray*}
u(e_{1,2},z,\vec x)&=&\Gamma\left(\frac {z-x_2}{x_3},
\frac{x_1-x_2}{x_3},-\frac{x_1}{x_3}\right)^{-1},\\
u(e_{3,2},z,\vec x)&=&\Gamma\left(\frac {z}{x_1},\frac{x_2-x_3}{x_1},
\frac{x_3}{x_1}\right),\\
u(e_{i,j},z,\vec x)&=&1,\qquad j\neq 2,\\
u(t_2,z,\vec x)&=&\theta_0\left(\frac{z-x_2}{x_1},\frac{x_3}{x_1}\right),\\
u(t_j,z,\vec x)&=&1,\qquad j\neq2,
\end{eqnarray*}
extend to a 1-cocycle $u:G\to N/M$. Its cohomology class 
$[u]\in H^1(G,N/M)$
is independent of the choice of extension and is  non-trivial.
  The corresponding
cohomology class $[\phi]=\delta_*[u]$ is represented by a function
$\phi:G^2\to M$ whose restriction to $(\Z^3)^2$ is given by
\begin{equation}\label{e-deltastar}
\phi(t_1^lt_2^mt_3^n,t_1^{l'}t_2^{m'}t_3^{n'},z,\vec x)
=
e^{\pi i(nm'(2z/x_1+1)-m'n(n+1)x_3/x_1-nm'(m'+1+2m)x_2/x_1)}.
\end{equation}
\end{thm}
\medskip

\noindent{\it Proof:} 
The proof is based a the presentation of $G$ by generators
and relations. The $\mathrm{SL}(3,\Z)$ subgroup is generated
by the elementary matrices $e_{i,j}$, ($i\neq j$). The relations
can be chosen  \cite{M} to be
\begin{eqnarray*}
e_{i,j}e_{k,l}&=&e_{k,l}e_{i,j},\qquad i\neq k,\quad j\neq l,\\
e_{i,j}e_{j,k}&=&e_{i,k}e_{jk}e_{i,j},\\
(e_{1,3}e_{3,1}^{-1}e_{1,3})^4&=&1.
\end{eqnarray*}
The relations of the generators of the $\Z^3$ subgroup are
$t_it_j=t_jt_i$ and the relations between $e_{i,j},t_k$ are
\begin{eqnarray*}
e_{i,j}t_k&=&t_ke_{i,j},\qquad i\neq k,\\
e_{i,j}t_i&=&t_it_j^{-1}e_{i,j}.
\end{eqnarray*}
The cocycle condition uniquely determines a 1-cocycle in terms
of its values on generators. For any functions
$u_{i,j}, u_k\in N/M$ there exists 
a unique 1-cocycle of the free group on generators
$e_{i,j},e_{k}$ such that $u_{i,j}=u(e_{i,j})$ and $u_i=u(t_i)$.
 This cocycle defines a 1-cocycle of $G$ if and only if
the relations are sent to 1. 

This can be checked using the identities of the functions
$\Gamma$, $\theta_0$. Let us consider some non-trivial examples.

The relation $t_2e_{3,2}t_3=t_3e_{3,2}$ translates to the
condition
\[
u_2(z,\vec x)u_{3,2}(z-x_2,\vec x)u_3(z-x_2,e_{3,2}^{-1}\vec x)
=u_3(z,\vec x)u_{3,2}(z-x_3,\vec x),
\]
which reduces to the defining functional equation for $\Gamma$:
\[
\theta_0
\left(
\frac{z-x_2}{x_1},
\frac{x_3}{x_1}
\right)
\Gamma
\left(
\frac{z-x_2}{x_1},
\frac{x_2-x_3}{x_1},
\frac{x_3}{x_1}
\right)
=\Gamma
\left(
\frac{z-x_3}{x_1},
\frac{x_2-x_3}{x_1},
\frac{x_3}{x_1}
\right).
\]
The relation $e_{1,3}e_{3,2}=e_{1,2}e_{3,2}e_{1,3}$ translates
to the condition
\[
u_{1,3}(z,\vec x)u_{3,2}(z,e_{1,3}^{-1}\vec x)
=
u_{1,2}(z,\vec x)u_{3,2}(z,e_{1,2}^{-1}\vec x)u_{1,3}(z,e_{1,2}^{-1}
e_{1,3}^{-1}\vec x)\mod M.
\]
By inserting the given expressions for $u_{i,j}$, we see
that the condition is
\begin{equation}\label{e-432}
\Gamma
\left(
\frac{z}{x_1-x_3},
\frac{x_2-x_3}{x_1-x_3},
\frac{x_3}{x_1-x_3}
\right)
=
\frac{
\Gamma
\left(
\frac{z}{x_1-x_2},
\frac{x_2-x_3}{x_1-x_2},
\frac{x_3}{x_1-x_2}
\right)}
{\Gamma
\left(
\frac{z-x_2}{x_3},
\frac{x_1-x_2}{x_3},
-\frac{x_1}{x_3}
\right)}\quad \mod M
.
\end{equation}
Using the fact that $\Gamma$ is periodic with period 1 in
all its arguments, and setting $Z=(z-x_1+x_3)/(x_1-x_3)$,
$\sigma=x_3/(x_1-x_3)$, $\tau=(x_2-x_1)/(x_1-x_3)$,
one sees that this identity reduces to
\[
\Gamma
\left(
{Z},
\tau,
\sigma
\right)
=
\frac{
\Gamma
\left(
-\frac{Z+1}{\tau},
-\frac{1}{\tau},
-\frac{\sigma}{\tau}
\right)}
{\Gamma
\left(
\frac{Z-\tau}{\sigma},
-\frac{\tau}{\sigma},
-\frac{1}{\sigma}
\right)}\quad \mod M.
\]
By the last identity of Proposition \ref{p-3} and the
rules for changing signs before Theorem \ref{t-3}, this identity
reduces to one of 
the three-term relation in Theorem \ref{t-4}.

To show that the class of $[\phi]=\delta_*[u]$ is non-trivial, one
notices that the restriction of $[\phi]$ to $\Z^3$ is
given by the formula \Ref{e-deltastar}. So its image under
the map $c_1:H^1(\Z^3,M)\to H^3(\Z^3,2\pi i\Z)$, which comes
from the short exact sequence \[0\to2\pi i\Z\to 2\pi i\Q(x_1/x_3,x_2/x_3)[z/x_3]
\stackrel\exp\to M\to 1,\] is calculated as in the proof 
of Prop.\ \ref{p-sic} and is non-trivial.
$\square$

\subsection{Explicit description of the 2-cocycle}

To describe the cohomology class $[\phi]=\delta_*[u]\in H^2(G,M)$ 
of $G=\mathrm{SL}(3,\Z)\semidirect\Z^3$ with coefficients $M$
 arising in Theorem \ref{t-co},
it is convenient to use the isomorphism of $H^2(G,M)$ with the set 
$\mathcal{E}(G,M)$ of
equivalence classes of group extensions
\[
1\to M\stackrel{i}\to E\stackrel{p}\to G\to 1,
\]
of the $G$-module $M$. The isomorphism assigns to any such extension
its {\em characteristic class} in $H^2(G,M)$. It is defined as follows:
choose a map $\sigma:G\to E$ such that $p\circ \sigma=\mathrm{id_G}$. Then
$\sigma(g)\sigma(h)=i(\phi(g,h))\sigma(gh)$ for some 2-cocycle 
$\phi\in C^2(G,M)$ whose class, the characteristic class of the extension,
is independent of the choice of $\sigma$.

In our case the extension can be described in terms of generators
and relations.

Let us introduce a set of elements 
$\phi_{j}^k=\exp(i\pi L_{j}^k)$,
$\phi_{j,k}^l=\exp(i\pi L_{j,k}^l)$,
$\phi_{j,k}^{l,m}=\exp(i\pi L_{j,k}^{l,m})$, ($1\leq j,k,l,m\leq 3$)
of $M$ with
\[
L_{1,2}^{2}= {\relax 
\frac {x_2(6z^{2} -6( {x_{3}} + 2{x_{2}})
z + {x_{2}}{x_{1}} - {x_{1}}^{2} + 6{x_{2}}^{2}
 + 6{x_{2}}{x_{3}} + {x_{3}}^{2})}{6{x_{3}}{x_{1}}({x_{2}}
 - {x_{1}})}} \,,
\]
\[
L_{1,2}^{1}={ 
\frac { - 6z^{2} + 6({x_{1}} + {x_{3}} + 2{x_{2}})z
 - {x_{3}}^{2} - 6{x_{2}}{x_{1}} + 3{x_{1}}{x_{3}} - 6{
x_{2}}{x_{3}} - 6{x_{2}}^{2} - {x_{1}}^{2}}{6{x_{1}}{x_{3}}}
} \,,
\]
\[
L_{1,3}^{2}= {\relax 
\frac {6z^{2} -6 ( {x_{3}} + 2{x_{2}})z + 6{x_{2}}
{x_{3}} + 5{x_{1}}^{2} + 6{x_{2}}^{2} + {x_{3}}^{2} - 5{x_{1}}{x_{3}}}
{6({x_{3}} - {x_{1}}){x_{1}}}} \,,
\]
\[
L_{2}^{3}={\relax \frac {2z - 2{x_{2}} - 2{x_{3}}
 + {x_{1}}}{{x_{1}}}} \,,
\]
\[
L_{1,3}^{3,2}=
 {\relax 
\frac {( 2z- {x_{2}})(2z^{2} - 2z{x_{2}} - {x_{1}}^{
2} + {x_{1}}{x_{3}} - {x_{3}}^{2} + {x_{2}}{x_{1}})}
{12({x_{3}} - {x_{1}}){x_{3}}({x_{2}} - {x_{1}})}} \,,
\]
\[
L_{3,1}^{1,2}={\relax} {\relax 
\frac {( 2z- {x_{2}} )(2z^{2} - 2z{x_{2}} + {x_{2}}
{x_{3}} - {x_{1}}^{2} + {x_{1}}{x_{3}} - {x_{3}}^{2})}{12{x_{1}}({x_{2}} - {x_{3}})({x_{3}} - {
x_{1}})}} \,,
\]
\[
L_{1,2}^{3,2}={\relax} {\relax 
\frac {{x_{2}}(2z - {x_{2}})(2z^{2} - 2z{x_{2}} + 
{x_{2}}{x_{3}} - {x_{3}}^{2} - {x_{1}}^{2} + {x_{2}}{x_{1}})}{12{x_{1}}({x_{2}} - {x_{1}}){x_{3}}
({x_{2}} - {x_{3}
})}} \,,
\]
\[
L_{3,2}^{1,2}=-L_{1,2}^{3,2},\qquad L_{3}^2=-L_{2}^3,
\]
and all other $L_{\cdots}^{\cdots}$ are zero.
Let $E$ be the group generated by  $\hat e_{i,j}$, $1\leq i\neq j\leq3$, 
$\hat t_i$,
$1\leq i\leq3$ and the elements of  $M$ subject to the following defining relations:

\noindent
a) The product of elements of $M$ in $E$ is the product in $M$.

\noindent
b) Relations with $M$:
\[
\hat e_{i,j}u=\rho(e_{i,j})(u) \hat e_{i,j}, \quad \hat t_iu=\rho(t_i)(u)\hat t_i,
\] 
$u\in M$.

\noindent c) Relations among the $\hat e_{i,j}$:
\begin{eqnarray*}
\hat e_{i,j}\hat e_{k,l}&=&\phi_{i,j}^{k,l}\hat e_{k,l}\hat e_{i,j},\qquad i\neq k,\quad j\neq l,\\
\hat e_{i,j}\hat e_{j,k}&=&\phi_{i,j}^{j,k}\hat e_{i,k}\hat e_{j,k}\hat e_{i,j},\\
(\hat e_{1,3}\hat e_{3,1}^{-1}\hat e_{1,3})^4&=&1.
\end{eqnarray*}
d) Relations among the $\hat t_i$: $\hat t_i\hat t_j=\phi_{i}^{j}\hat t_j\hat t_i$
\par\noindent\noindent
e) Relations between $\hat e_{i.j}$ and $\hat t_i$:
\begin{eqnarray*}
\hat e_{i,j}\hat t_k&=&\phi_{i,j}^k\hat t_k\hat e_{i,j},\qquad i\neq k,\\
\hat t_j\hat e_{i,j}\hat t_i&=&\phi_{i,j}^i\hat t_i\hat e_{i,j}.
\end{eqnarray*}
This group comes with a natural homomorphism $M\to E$.

\begin{thm} Let $G, M, [\phi]$ be as in Theorem \ref{t-co}.
The map $M\rightarrow E$ fits into a group extension
\[
1\to M\to E\to G\to 1,
\]
of the $G$-module $M$,
whose characteristic class in $H^2(G,M)$ 
is  $[\phi]=\delta_*[u]$ of Theorem
\ref{t-co}.
\end{thm}

\medskip
To prove this theorem we need to recall some facts about the
description of extensions by generators and relations. Let
$G=F/R$ with $R$ a normal subgroup of a free group $F$ 
with generators $(e_i)_{i\in I}$.  The canonical projection
$F\to G$ will  be denoted by $x\mapsto \bar x$.
Let $M$ be a $G$-module and
$\rho:G\to \mathrm{Aut}(M)$ the corresponding homomorphism.

Suppose that a cohomology class $[\phi]\in H^2(G,M)$ is given.
We want to describe the middle group of the corresponding 
extension by generators
and relations. The relations are
written in terms of a map $\psi:R\to M$ built out of $[\phi]$. 
We proceed to explain
how to construct an extension associated to a map
 $\psi$ with certain properties, and how to construct 
$\psi$ given the characteristic class $[\phi]$ of the
extension.

Abstractly, the relation between $\psi$ and $[\phi]$ 
is that $\psi$ is any inverse image of $[\phi]$ by
the surjective homomorphism
$H^0(G,\mathrm{Hom}(R_{\mathrm{ab}},M))\to H^2(G,M)$
described in \cite{MacL}, Section VIII.9. 
Here $R_{\mathrm{ab}}$ denotes the abelianization $R/[R,R]$ of $R$,

This
homomorphism can be described explicitly as follows. 

An element of 
$H^0(G,\mathrm{Hom}(R_{\mathrm{ab}},M))$
is, by definition, a map $\psi:R\to M$ such that 
\begin{equation}\label{e-rr}
\psi(rs)=\psi(r)\psi(s)\quad \mathrm{and}\quad \psi(xrx^{-1})
=\rho(\bar x)\psi(r),
\end{equation}
for all $r,s\in R$, $x\in F$. The semidirect
product $\hat E=F\times_{\rho} M$ is the cartesian product with group
multiplication rule $(x,u)(y,v)=(xy,u\rho(x)v)$. Then the properties of
$\psi$ imply that $\hat R=\{(r,\psi(r)^{-1})\,|\, r\in R\}$ is a normal
subgroup of $\hat E$. Let $E=\hat E/\hat R$. Then it is easy to show
that we have an extension
$1\to M\to E\to G\to 1$, with the obvious maps. If $R$ is generated by
relations $r_j$, $j\in J$,
 then $E$ is the group with generators $(e_i)_{i\in I}$, $M$ and defining
relations: 
\begin{enumerate}\item the product of elements in $M$ in $E$ is the product in $M$,
\item $e_iu=(\rho(e_i)u)e_i$,
\item $r_j=\psi(r_j)$.
\end{enumerate}
$i\in I, j\in J,u,v\in M$.
The characteristic class of this extension is the class of $\psi\circ \mu$,
where $\mu(g,h)$ is defined by $\lambda(g)\lambda(h)=\mu(g,h)\lambda(gh)$,
for any section $\lambda:G\to F$, as can be seen by choosing the section
$\sigma(g)=[(\lambda(g),1)]\in E$.

Conversely, given a 2-cocycle $\phi\in C^2(G,M)$  we may find 
a $\psi\in H^0(G,\mathrm{Hom}(R_{\mathrm{ab}},M))$ mapping to 
$[\phi]$ as follows. Let $\psi:F\to M$ be the unique map
such that 
\begin{eqnarray}\label{e-99b}
\psi(e_i^{\pm1})&=&1, \qquad i=1,\dots,n,\\
\psi(xy)&=&\phi(\bar x,\bar y)\psi(x)\rho(\bar x)\psi(y), \qquad \forall x,y\in F.\notag
\end{eqnarray}
Then the restriction of $\psi$ to $R$ obeys \Ref{e-rr}. For
any section $\lambda:G\to F$ we have
\begin{eqnarray*}
\psi(\mu(g,h))\psi(\lambda(gh))&=&\psi(\mu(g,h)\lambda(gh))\\
 &=&\psi(\lambda(g)\lambda(h))\\
 &=&\phi(g,h)\psi(\lambda(g))\rho(g)\psi(\lambda(h)).
\end{eqnarray*}
Thus the 2-cocycle $\psi\circ \mu$
is indeed in the same cohomology class as $\phi$.

Let us apply this construction to our case. Let $[u]\in H^1(G,N/M)$ be
the cohomology class described in the Theorem. Let us choose a
representative $u:G\to N$, which on the generators 
coincides with
the functions given in the claim. Then a representative of 
the class $[\phi]=\delta_*[u]\in H^2(G,M)$
is given by $\phi(g,h)=u(gh)/(u(g)\rho(g)u(h))$.
The generators
$e_i$ are here $e_{i,j},t_i$. We need to compute the
value of $\psi$ on relations.

Every element $\neq 1$ of $F$ is can uniquely be 
written as a reduced word $x_1\cdots x_k$, with $x_j\in\{e_i^{\pm1},\,i\in I\}$.
Reduced means that $x_j\neq x_{j+1}^{-1}$, $j=1,\dots k-1$. The function
$\psi:F\to M$ is then given according to the rule \Ref{e-99b}
by the formula
\[
\psi(x_1\cdots x_k)=\frac{u(\bar x_1)\prod_{i=2}^k\rho(\bar x_1\cdots 
\bar x_{i-1})u(\bar x_i)}
{u(\bar x_1\cdots \bar x_k)}\,,
\]
for $x_i\in\{e_j,e_j^{-1}\}$. If $x_1\cdots x_k$ is a relation,
the denominator is equal to 1, so to compute $\psi(r)$ we only need
the value of $u$ on generators and their inverses. 
The value of $u$ on inverses of generators obeys, by the cocycle condition, 
$u(\bar x^{-1})=(\rho(\bar x^{-1})u(\bar x))^{-1}$ modulo $M$.
Since we have no generators such that $\bar x=\bar x^{-1}$ we
may choose $u$ so that this relation holds for representatives,
not just modulo $M$.

We explain the calculation in the case of the relation 
$r=e_{1,3}e_{3,2}(e_{1,2}e_{3,2}e_{1,3})^{-1}$. For the other
relations the calculation is similar.
The equation for $\phi_{1,3}^{3,2}=\psi(r)$ can be written as
\[
u_{1,3}(z,\vec x)u_{3,2}(z,e_{1,3}^{-1}\vec x)
=
\psi(r,z,\vec x)
u_{1,2}(z,\vec x)u_{3,2}(z,e_{1,2}^{-1}\vec x)u_{1,3}(z,e_{1,2}^{-1}
e_{1,3}^{-1}\vec x).
\]
with $u_{i,j}=u(e_{i,j})$. 
\[
\Gamma
\left(
\frac{z}{x_1-x_3},
\frac{x_2-x_3}{x_1-x_3},
\frac{x_3}{x_1-x_3}
\right)
=
\psi(r,z,\vec x)\,
\frac{
\Gamma
\left(
\frac{z}{x_1-x_2},
\frac{x_2-x_3}{x_1-x_2},
\frac{x_3}{x_1-x_2}
\right)}
{\Gamma
\left(
\frac{z-x_2}{x_3},
\frac{x_1-x_2}{x_3},
-\frac{x_1}{x_3}
\right)}
.
\]
Comparing with Theorem \ref{t-4}, we see that (cf.\  the calculation
following \Ref{e-432})
\[
\psi(r,z,\vec x)=\exp\left(-\pi i 
Q\left(
\frac{z-x_1+x_3}{x_1-x_3};\frac{x_2-x_1}{x_1-x_3},
\frac{x_3}{x_1-x_3}\right)\right).
\]
This expression is easily checked to be identical
to $\exp(\pi i L_{1,3}^{3,2})$.

Proceeding in the same way with the other relations, we find
that the non-trivial values of $\phi$ on relations 
are $\phi_A^B=\exp(\pi iL_A^B)$ with
\[
{L_{1,2}^{2}} =  - 2\,{\displaystyle \frac {z - {x_{2}}}{{x
_{3}}}}  + 1 + {F}\left({\displaystyle \frac {z - {x_{2}}}{{x_{
1}}}} , \,{\displaystyle \frac {{x_{3}}}{{x_{1}}}} \right) - {F}
\left({\displaystyle \frac {z - {x_{1}}}{{x_{1}} - {x_{2}}}} , \,
{\displaystyle \frac {{x_{3}}}{{x_{1}} - {x_{2}}}} \right),
\]
\[
{L_{1, 2}^{1}} = 2\,{\displaystyle \frac {z - {x_{2}}}{{x_{3}
}}}  + 1 - {F}\left({\displaystyle \frac {z - {x_{2}}}{{x_{1}}}} , \,
{\displaystyle \frac {{x_{3}}}{{x_{1}}}} \right),\qquad
{L_{1, 3}^{2}} =  - {F}\left({\displaystyle \frac {z - {x_{
2}}}{{x_{1}} - {x_{3}}}} , \,{\displaystyle \frac {{x_{1}}}{{x_{1
}} - {x_{3}}}} \right),
\]
\[
{L_{2}^{3}} = 2\,{\displaystyle \frac {z - {x_{2}} - {x_{3}}}{{
x_{1}}}}  + 1,\qquad
{L_{3,1}^{1,2}} = {Q}\left({\displaystyle \frac {z - {x
_{1}}}{{x_{1}}}} ; \,{\displaystyle \frac 
{{x_{2}} - {x_{3}}}{{x_{1}}}} , \,{\displaystyle \frac {{x_{3}} - {x_{1}}}{{x_{1}}}} \right),
\]
\[
{L_{1, 2}^{3, 2}} = {Q}\left({\displaystyle \frac {z - {x
_{1}}}{{x_{1}}}} ; \,{\displaystyle \frac {{x_{2}} - {x_{3}}}{{x_{1}}}} , 
\,{\displaystyle \frac {{x_{3}} - {x_{1}}}{{x_{1}}}} \right) + {Q}\left(
{\displaystyle \frac {z - {x_{1}} + {x_{3}}}{{x_{1}} - {x_{3}}}} ;
 \,{\displaystyle \frac {{x_{3}}}{{x_{1}} - {x_{3}}}} , \,{\displaystyle \frac {{x_{2}} - {x_{1}}}
{{x_{1}} - {x_{3}}}} \right).
\]
Here $F$ is the polynomial appearing in the modular transformation
properties of $\theta_0$:
\[
{F}\left(z, \,\tau \right)={\displaystyle \frac {z^{2}}{\tau }}  + z
\,\left({\displaystyle \frac {1}{\tau }}  - 1\right) + {\displaystyle 
\frac {\tau}{6}}  + {\displaystyle \frac {1}{2}}  + 
{\displaystyle \frac {1}{6\tau }} 
\]
These expressions may be (preferably with a computer) simplified
to give the claim of the theorem. The proof is complete.

\subsection{The restriction of the 2-cocycle to $\mathrm{SL}(3,\Z)$}
Let $G,M,[\phi]$ be as in Theorem \ref{t-co}.
We show here that the restriction $[\bar\phi]\in H^2(\mathrm{SL}(3,\Z),M)$
to $\mathrm{SL}(3,\Z)\subset G$  is non-trivial
(See \cite{S} for a description of the cohomology of $\mathrm{SL}(3,\Z)$).
This is proved by showing that the restriction to
a $D_4$ subgroup is non-trivial. This $D_4$ subgroup is generated by
$a=(e_{2,1}^{-1}e_{1,2}e_{2,1}^{-1})^2:(x_1,x_2,x_3)\mapsto
(-x_1,-x_2,x_3)$
and
$b=e_{3,1}e_{1,3}^{-1}e_{1,3}:(x_1,x_2,x_3)\mapsto (-x_3,x_2,x_1)$
The defining relations of $D_4$ are 
$ a^2= b^4=1$, $b a b= a $.
The restriction of $[\phi]$ to this subgroup is the characteristic
class of the extension $1\to M\to\hat D_4\to D_4$ with
$\hat D_4$ is the inverse image of $D_4$ by the projection
$E\to G$. A presentation of $\hat D_4$ is obtained by
choosing lifts of generators
\begin{eqnarray*}
\hat a&=&e^{\frac{\pi i}2\left(\frac {z^2}{x_1x_3}-\frac{2z}{x_3}+1
+\frac{x_1}{6x_3}+\frac{x_3}{6x_1}\right)}(\hat e_{2,1}^{-1}
\hat e_{1,2}\hat e_{2,1}^{-1})^2
\\
\hat b&=&\hat e_{3,1}\hat e_{1,3}^{-1}\hat e_{1,3}.
\end{eqnarray*}
The relations between these generators  are computable
from the presentation of $E$ above, with the result
\[
\hat a^2=\hat b^4=1,\qquad \hat b\hat a\hat b=i\hat a,\qquad i=e^{\frac{2\pi i}4}\in M.
\]
\begin{proposition} The pull-back $i^*[\phi]$ of the class $[\phi]\in H^2(G,M)$
of Theorem \ref{t-co} by the inclusion $i:D_4\to SL(3,\Z)\semidirect\Z^3$ 
is the characteristic class in $H^2(D_4,M)$ of the extension
$1\to M\to \hat D_4\to D_4$. It is a non-trivial cohomology class.
\end{proposition}

\medskip
\noindent
{\it Proof.} 
It remains to prove that the class is non-trivial.
If the characteristic class were trivial, the exact sequence
of the extension would split. This would mean that, for suitable
homogeneous functions $A,B\in\Q(x_1,x_2,x_3)[z]$,
$e^{2\pi iA}\hat a$ and $e^{2\pi i B}\hat b$ obey the relations
of $D_4$. 
Suppose, by contradiction, that such functions exist.
Then $A,B$  obey
\begin{gather*}
A(z,\vec x)+A(z,a^{-1}\vec x)=r,
\\
B(z,\vec x)+B(z,b^{-1}\vec x)+B(z,b^{-2}\vec x)+B(z,b^{-3}\vec x)=s,
\\
B(z,\vec x)+A(z,b^{-1}\vec x)+B(z,(ba)^{-1}\vec x)=-\frac14+A(z,\vec x)
+t.
\end{gather*}
for some integers $r,s,t$. Let $\bar A=A-r/2$, $\bar B=B-s/4$, then
$\bar A$, $\bar B$ obey
\begin{gather*}
{\bar A}(z,\vec x)+{\bar A}(z,a^{-1}\vec x)=0,
\\
{\bar B}(z,\vec x)+{\bar B}(z,b^{-1}\vec x)+{\bar B}(z,b^{-2}\vec x)+{\bar B}(z,b^{-3}\vec x)=0,
\\
{\bar B}(z,\vec x)+{\bar A}(z,b^{-1}\vec x)+{\bar B}(z,(ba)^{-1}\vec x)=-\frac14+{\bar A}(z,\vec x)
+t-\frac s2.
\end{gather*}
If we view $\Q(x_1,x_2,x_3)[z]$ as a module over the group ring
$\Z D_4$, the first two equations can be written as $(1+a){\bar A}=0$,
$(1+b+b^2+b^3){\bar B}=0$. This implies that ${\bar A}$, ${\bar B}$ are annihilated 
by the idempotent
\[
P=\frac18\sum_{g\in D_4}g
=\frac18(1+b+b^2+b^3)(1+a)
=\frac18(1+a)(1+b+b^2+b^3).
\]
Applying $P$ to the third equation, we get $0=-1/4+t-s/2$, $t,s\in\Z$, a contradiction. $\square$

\end{document}